\documentclass{article}
\usepackage[leqno]{amsmath}
\usepackage{pstricks}
\usepackage[dvipdfmx]{graphicx}
\usepackage[english]{babel}
\usepackage{amsmath,enumerate, amsfonts, amssymb,amsthm, xypic}

\DeclareMathOperator{\R}{\mathbf R}
\DeclareMathOperator{\C}{\mathbf C}
\DeclareMathOperator{\N}{\mathbf N}
\DeclareMathOperator{\Z}{\mathbf Z}

\DeclareMathOperator{\Int}{Int}

\DeclareMathOperator{\id}{id}

\DeclareMathOperator{\bdry}{bdry}

\DeclareMathOperator{\Reg}{Reg}

\DeclareMathOperator{\Ima}{Im}

\begin{document}
\if0\ \vskip9mm

\begin{figure}[!h]
\begin{pspicture}(-0.85,1.3)(10.5,5)
\rput(3,3.62){\tiny$x_0$}
\pscircle*(3,3.5){0.03}
\rput(4.5,3.3){\tiny$2x_0$}
\rput(4.55,2.9){\tiny slope}
\rput(4.55,2.7){\tiny $=\!-1$}
\pscircle*(4,3.5){0.03}
\rput(4.6,3.9){\tiny$3x_0/2$}
\pscircle*(3.48,3.5){0.03}
\psline[linestyle=dashed, dash=2pt 1pt](2,2)(2,5)
\pscustom[linewidth=0.1pt]{\psline[linestyle=dashed, dash=0.1pt 0.1pt](2,2)(2,5)\gsave
\psline[linewidth=1pt](3,5)(3,4.5)(2.5,4)(2.5,3)(3,2.5)(3,2)
\fill[linewidth=0.1pt,fillstyle=hlines]\grestore}
\psline(2,5)(3,5)
\psline(2,2)(3,2)
\psline(3,5)(3,4.5)(2.5,4)(2.5,3)(3,2.5)(3,2)
\psline(3,4.5)(3.48,4.02)(3.48,2.98)(3,2.5)
\psline(2.7,3.8)(2.7,3.2)(3,2.9)(3.3,3.2)(3.3,3.8)(3,4.1)(2.7,3.8)
\psline[linestyle=dotted, dash=1pt 0.1pt](2.5,3)(4.15,4.65)
\psline[linestyle=dotted, dash=1pt 0.1pt](2.5,4)(4.15,2.35)
\psline[linestyle=dotted, dash=1pt 0.1pt](3,2.5)(4.15,2.5)
\psline[linestyle=dotted, dash=1pt 0.1pt](3,4.5)(4.15,4.5)
\psline(3.2,5)(3.2,4.5)(3.62,4.12)(3.62,2.88)(3.2,2.5)(3.2,2)
\psline(3.4,5)(3.4,4.5)(3.74,4.2)(3.74,2.78)(3.4,2.5)(3.4,2)
\rput(3.1,1.3){\footnotesize$x_0\in(0,\,2/3)\!\times\!\{0\}$}
\psline(4,5)(4,2)
\psline(4.15,5)(4.15,2)
\rput(1.9,1.65){\scriptsize$\{0\}\!\times\![-1,1]$}
\rput(4.25,1.65){\scriptsize$\{1\}\!\times\![-1,1]$}
\rput(1.6,3.93){\scriptsize{value}}
\rput(1.55,3.7){\scriptsize=}
\rput(1.58,3.5){\scriptsize$|x_0|/2$}
\rput(1.6,2.6){\tiny{slope}}
\rput(1.6,2.4){\tiny$=1$}
\pscurve[linestyle=dashed, dash=1pt 1pt]{->}(1.7,3.3)(1.85,3.2)(2.25,3.2)
\pscurve[linestyle=dashed, dash=1pt 1pt]{->}(4.45,3.43)(4.25,3.5)(4.05,3.5)
\psline[linestyle=dashed, dash=1pt 1pt]{->}(4.5,3.8)(3.55,3.5)
\pscurve[linestyle=dashed, dash=1pt 1pt]{->}(1.9,2.5)(2.5,2.65)(2.65,3.1)

\pscurve[linestyle=dashed, dash=1pt 1pt]{->}(4.5,2.8)(4.25,2.85)(4.05,2.7)
\psline[linestyle=dashed, dash=2pt 1pt](6.2,5)(6.2,2)
\pscustom[linewidth=0.1pt]{\psline[linestyle=dashed, dash=0.1pt 0.1pt](6.2,2)(6.2,5)\gsave
\psline[linewidth=1pt](7.7,5)(7.7,4.6)(7.05,3.95)(7.05,3.04)(7.7,2.4)(7.7,2)
\fill[linewidth=0.1pt,fillstyle=hlines]\grestore}
\psline(6.2,5)(7.7,5)
\psline(6.2,2)(7.7,2)
\psline(7.7,4.6)(8.17,4.13)
\psline(8.17,2.87)(7.7,2.4)
\rput(7.7,3.62){\footnotesize$x_0$}
\psline(7.7,5)(7.7,4.6)(7.05,3.95)(7.05,3.04)(7.7,2.4)(7.7,2)
\rput(7.25,1.3){\footnotesize$x_0\in[2/3,\,1)\!\times\!\{0\}$}
\pscircle*(7.7,3.5){0.04}
\psline(7.3,3.9)(7.3,3.1)(7.7,2.7)(8.1,3.1)(8.1,3.9)(7.7,4.3)(7.3,3.9)
\psline[linestyle=dotted, dash=1pt 0.1pt](8.17,3.97)(7.7,3.5)
\psline[linestyle=dotted, dash=1pt 0.1pt](7.7,3.5)(8.17,3.03)
\psline(7.9,5)(7.9,4.6)(8.17,4.33)
\psline(7.9,2)(7.9,2.4)(8.17,2.67)
\psline(8.1,5)(8.1,4.6)(8.17,4.53)
\psline(8.1,2)(8.1,2.4)(8.17,2.47)
\rput(6.2,1.65){\scriptsize$\{0\}\!\times\![-1,1]$}
\rput(8.17,1.65){\scriptsize$\{1\}\!\times\![-1,1]$}
\rput(5.8,3.93){\scriptsize{value}}
\rput(5.75,3.7){\scriptsize=}
\rput(5.77,3.5){\scriptsize$|x_0|/2$}
\pscurve[linestyle=dashed, dash=1pt 1pt]{->}(5.9,3.3)(6.1,3.2)(6.5,3.2)

\psline[linestyle=dashed, dash=2pt 1pt](10,5)(10,2)
\pscustom[linewidth=0.1pt]{\psline[linestyle=dashed, dash=0.1pt 0.1pt](10,2)(10,5)\gsave
\psline[linewidth=1pt](12,5)(12,4.5)(11,4)(11,3)(12,2.5)(12,2)
\fill[linewidth=0.1pt,fillstyle=hlines]\grestore}
\psline(10,5)(12,5)
\psline(10,2)(12,2)
\psline(12,5)(12,4.5)(11,4)(11,3)(12,2.5)(12,2)
\psline(12,4.15)(11.35,3.8)(11.35,3.15)(12,2.85)
\psline(12,3.85)(11.7,3.7)(11.7,3.3)(12,3.15)
\rput(9.9,1.65){\scriptsize$\{0\}\!\times\![-1,1]$}
\rput(12.1,1.65){\scriptsize$\{1\}\!\times\![-1,1]$}
\rput(12,3.62){\footnotesize$x_0$}
\rput(11.1,1.3){\footnotesize$x_0=(1,0)$}
\pscircle*(12,3.5){0.04}
\rput(9.6,3.9){\scriptsize{value}}
\rput(9.6,3.6){\scriptsize{=$1/2$}}
\pscurve[linestyle=dashed, dash=1pt 1pt]{->}(9.7,3.4)(9.85,3.3)(10.4,3.3)

\end{pspicture}
\caption{Level sets of the function $x\to d_X(x,x_0)$}
\begin{quote}{\small Each level curve of $x\to d_X(x,x_0)$ is a union of vertical segments and segments with slope $\pm1$ if the value is not $|x_0|/2$.
If $d_X(x,x_0)=x_0/2$, each level set is a union of such three segments and a two-dimensional domain.}\end{quote}
\end{figure}

\begin{figure}[!h]
\begin{pspicture}(-0.55,1.6)(10.5,5)
\rput(3.4,3.5){\small0}
\pscircle(3.4,3.5){1.5}
\pscircle(3.4,3.5){1}
\pscircle(3.4,3.5){0.5}
\rput(3.4,1.6){$b_0=0$}
\rput(4.8,5){\scriptsize{value}}
\rput(4.8,4.7){\scriptsize{=$1/2$}}
\psline[linestyle=dashed, dash=1pt 1pt]{->}(4.8,4.54)(4.65,4.34)

\rput(7.2,3.5){\small0}
\pscircle(7.2,3.5){1.5}
\rput(7.2,1.6){$|b_0|<2/3$}
\pscustom[linewidth=0.7pt,fillstyle=hlines]{\psarc(7.2,3.5){0.35}{5}{45}
\psarc(7.2,3.5){0.6}{50}{0}}
\psarc(7.2,3.5){0.9}{-4}{54}
\pscurve(7.55,3.55)(7.6875,3.55)(7.825,3.5)(7.925,3.5)(8.1,3.42)

\pscurve(7.42,3.76)(7.5025,3.83)(7.605,3.985)(7.67,4.05)(7.73,4.2)
\psarc(7.2,3.5){1.15}{-15}{65}
\psarc(7.2,3.5){0.9}{65}{-15}
\pscurve(7.58,4.3)(7.65,4.42)(7.66,4.55)
\pscurve(8.06,3.29)(8.22,3.25)(8.32,3.17)
\pscircle(7.2,3.5){1.35}
\pscircle*(7.77,3.77){0.03}
\rput(8.72,4.48){\small$b_0$}
\psarc(7.2,3.5){0.74}{15}{35}
\psarc(7.2,3.5){0.53}{15}{35}
\pscurve(7.71,3.63)(7.76,3.66)(7.81,3.665)(7.86,3.7)(7.91,3.7)
\pscurve(7.65,3.81)(7.6875,3.82)(7.725,3.87)(7.7625,3.88)(7.8,3.93)
\pscurve[linestyle=dashed, dash=1pt 1pt]{->}(8.55,4.4)(8,4)(7.8,3.8)

\rput(11,3.5){\small0}
\rput(11,1.6){$|b_0|=1$}
\pscircle(11,3.5){1.5}
\pscustom[linewidth=0.7pt,fillstyle=hlines]{\psarc(11,3.5){0.75}{-55}{102}
\pscurve(10.85,4.22)(10.64,4.55)(10.4,4.84)
\psarc(11,3.5){1.5}{115}{-67}
\pscurve(11.58,2.16)(11.54,2.505)(11.41,2.89)}
\psarc(11,3.5){1}{-25}{75}
\pscurve(11.3,4.97)(11.32,4.65)(11.26,4.43)
\pscurve(11.9,3.11)(12.15,2.87)(12.2,2.67)
\psarc(11,3.5){1.2}{0}{50}
\pscurve(11.75,4.4)(11.83,4.52)(11.85,4.71)
\pscurve(12.2,3.49)(12.35,3.48)(12.5,3.41)
\rput(12.57,4.15){\small$b_0$}
\pscircle*(12.36,4.1){0.04}
\end{pspicture}
\caption{Level sets of the function $b\to d(b,b_0)$}
\begin{quote}{\small Each level curve of $b\to d(b,b_0)$ is a union of curves of the form $\{c_1\exp(\pi i x_2):c_2\le x_2\le c_3\}$ and $\{(c_1\pm x_2)\exp(\pi i x_2):c_2\le x_2\le c_3\}$ if the value is not $|b_0|/2$.}\end{quote}
\end{figure}

\begin{figure}[!h]
\begin{pspicture}(-0.85,1.3)(10.5,5)
\rput(3,3.07){\footnotesize$x_0$}
\pscircle*(3,2.95){0.03}
\psline[linestyle=dashed, dash=2pt 1pt](2,2)(2,5)
\pscustom[linewidth=0.1pt]{\psline[linestyle=dashed, dash=0.1pt 0.1pt](2,2)(2,5)\gsave
\psline(3,5)(3,3.5)
\pscurve(3,3.5)(2.7,3.15)(2.5,2.8)
\psline(2.5,2.8)(2.5,2)
\fill[linewidth=0.1pt,fillstyle=hlines]\grestore}
\psline(2,2)(2.5,2)
\psline(2,5)(3,5)
\psline(2.65,2)(3,2)
\psline(3,5)(3,3.5)
\pscurve(3,3.5)(2.7,3.15)(2.5,2.8)
\psline(2.5,2.8)(2.5,2)
\pscurve(2.65,2)(2.75,2.23)(3,2.5)
\pscustom[linewidth=0.1pt]{\psline[linestyle=dashed, dash=2pt 1pt](2.65,2)(3,2.5)\gsave
\psline(3,2.5)(3,2)\fill[linewidth=0.1pt,fillstyle=hlines]\grestore}
\psline(3,2.5)(3,2)
\pscurve(3,3.5)(3.25,3.7)(3.5,3.8)
\psline(3.5,3.8)(3.5,2.8)
\pscurve(3.5,2.8)(3.25,2.7)(3,2.5)
\psline(3.4,5)(3.4,4.15)
\pscurve(3.4,4.15)(3.6,4.3)(3.9,4.4)
\psline(3.9,4.4)(3.9,2.65)
\pscurve(3.9,2.65)(3.6,2.49)(3.4,2.28)
\psline(3.4,2.28)(3.4,2)
\psline(4.2,5)(4.2,2)
\psline(4.05,5)(4.05,2)
\rput(1.9,1.65){\scriptsize$\{0\}\!\times\![-1,1]$}
\rput(4.25,1.65){\scriptsize$\{1\}\!\times\![-1,1]$}
\rput(3.1,1.3){\footnotesize$x_0\in(0,\,2/3)\!\times\!\{-1/3\}$}

\rput(7.2,3.5){\small0}
\pscircle(7.2,3.5){1.5}
\rput(7.2,1.6){$|b_0|<2/3$}
\pscustom[linewidth=0.7pt,fillstyle=hlines]{\psarc(7.2,3.5){0.35}{-60}{10}\gsave
\psarc(7.2,3.5){0.6}{30}{-13}\fill[linewidth=0.1pt,fillstyle=hlines]\grestore}
\psarc(7.2,3.5){0.6}{30}{-13}
\psarc(7.2,3.5){0.9}{-5}{43}
\pscurve(7.43,3.2)(7.62,3.31)(7.8,3.35)(7.9,3.43)(8.07,3.45)
\pscurve(7.54,3.52)(7.57,3.54)(7.6,3.6)(7.74,3.78)(7.86,4.13)
\psarc(7.2,3.5){1.15}{-21}{61}
\psarc(7.2,3.5){0.9}{55}{-25}
\pscurve(7.68,4.23)(7.75,4.35)(7.74,4.5)
\pscurve(8,3.11)(8.1,3.15)(8.3,3.1)
\pscircle(7.2,3.5){1.35}
\rput(8.7,4.46){\small$b_0$}
\pscircle*(7.79,3.57){0.03}
\pscurve[linestyle=dashed, dash=1pt 1pt]{->}(8.64,4.26)(8,3.62)(7.84,3.57)

\rput(11,3.5){\small0}
\rput(11,1.6){$|b_0|=1$}
\pscircle(11,3.5){1.5}
\pscustom[linewidth=0.7pt,fillstyle=hlines]{
\psarc(11,3.5){0.75}{-90}{90}\gsave
\pscurve(11,4.24)(10.6,4.4)(10.35,4.85)
\psarc(11,3.5){1.5}{115}{-65}\fill[linewidth=0.1pt,fillstyle=hlines]\grestore}
\pscurve(11,4.24)(10.6,4.4)(10.35,4.85)
\pscurve(11,2.73)(11.25,2.6)(11.6,2.15)
\psarc(11,3.5){1}{-40}{60}
\pscurve(11.43,4.92)(11.4,4.6)(11.5,4.37)
\pscurve(11.77,2.88)(11.97,2.95)(12.32,2.88)
\psarc(11,3.5){1.2}{-5}{45}
\pscurve(11.85,4.35)(11.85,4.45)(11.95,4.65)
\pscurve(12.2,3.38)(12.3,3.46)(12.5,3.48)
\rput(12.57,4.15){\small$b_0$}
\pscircle*(12.36,4.1){0.04}
\end{pspicture}
\caption{Level sets of $x\to\tilde d_X(x,x_0)$ and $b\to\tilde d(b,b_0)$}
\begin{quote}{\small Each level curve of $b\to\tilde d(b,b_0)$ is a union of curves of the form $\{c_1\exp(\pi i x_2):c_2\le x_2\le c_3\}$ and $\{x_1\exp\big(\pi i(c_4\pm x_1+\phi(x_1))\big):c_5\le x_2\le c_6\}$ if the value is not $|b_0|/2$.}\end{quote}
\end{figure}
\fi
\title{Semialgebraic Metric Spaces and Resolution of Singularities of Definable Sets}
\author{Masahiro Shiota}
\maketitle
\begin{abstract}
Consider the semialgebraic structure over the real field.
More generally, let an o-minimal structure be over a real closed field.
We show that a definable metric space $X$ with a definable metric $d$ is embedded into a Euclidean space so that its closure is compact and the metric on the image induced by $d$ is extended to a definable metric on the closure if and only if $\lim_{t,t'\to0}\allowbreak d(\gamma(t),\gamma(t'))=0$ for any definable continuous curve $\gamma:(0,\,1]\to X$ (Theorem 1).
We also find two compact semialgebraic metric spaces over the real field which are isometric but not semialgebraically isometric (Theorem 2).
A version of blow up is the key to the proof of Theorem 1.
Using it in the same way, we prove a resolution of singularities of definable sets (Theorem 3).
We prove the theorems by a constructive procedure.
\end{abstract}
\footnotetext{
\hspace{5.5mm}{\it Key words and phrases}. O-minimal structures, metric spaces, compactifications.}
\section{Introduction}
Research of general metric spaces is more difficult than that of compact metric spaces.
Hence we ask when a metric space is compactified in the sense that it is embedded in a compact metric space through a metric preserving map.
This problem is difficult in general, but it becomes simple in the o-minimal case.\par
Consider any o-minimal structure over any real closed field $R$.
Let a {\it definable metric space} mean a definable set with a definable (not necessarily Euclidean and always continuous) metric.
We call a definable metric space {\it compactifiable} if it can be embedded into a Euclidean space so that its closure is compact (i.e., bounded and closed in the Euclidean space) and the induced metric on the image is extendable to the closure.
We also call the closure or the embedding a {\it compactification}.
The simplest example of a non-compactifiable definable metric space is $R$ with the Euclidean metric.
We note that a definable non-compact manifold is homeomorphic (or diffeomorphic) to the interior of some compact definable manifold with boundary (Theorem II.5.7 in \cite{S1} and Theorem 2 in \cite{S3}), but a compactification of a definable noncompact metric manifold is not necessarily a manifold with boundary,
e.g., $R$ which is embedded in $R^2$ so that the image is $\{(x_1,x_2)\in R^2:|(x_1,x_2)|=1, (x_1,x_2)\not=(1,0)\}$ and whose metric is induced by the Euclidean metric on $R^2$, where the symbol $|\ \ |$ stands for the Euclidean norm.\vskip1mm\noindent
{\bf Theorem 1.} {\it A definable metric space $X$ with a metric $d$ is compactifiable if and only if $\lim_{t,t'\to0}\allowbreak d(\gamma(t),\gamma(t'))=0$ for any definable continuous curve $\gamma:(0,\,1]\to X$.
Moreover, a compactification is unique up to isometries.
To be precise, for any two compactifications $X_1$ and $X_2$, the identity map of $X$ induces an isometry between $X_1$ and $X_2$.}\vskip1mm
In the proof of Theorem 1 we will show that $d$ is bounded under the condition in Theorem 1.
However, there is an example of $X$ and $d$ where $d$ is bounded but the condition in Theorem 1 is not satisfied, e.g., $X=\R$ the real field, and $d(x,x')=\min\{|x-x'|,1\}$.
Namely, boundedness is not sufficient for compactification.
Even if we extend $\R$ to a real closed field $R$ in this example, the extended metric space over $R$ is not compactifiable.
An infinitude happens in this example, although the example is semialgebraic.
Another infinitude causes existence of two semialgebraic metric spaces which are isometric but not semialgebraically isometric (Theorem 2).\par
One of the simplest examples of compactifiable $(X,d)$, where $d$ cannot be extended to the closure $\overline X$, is the following:
Set $X=\{x\in\R^2:x\not=0,|x|\le 1\}$, let $\phi:Y\to\overline X$ be the blow up of $\overline X$ with center 0, let $Y\subset\R^n$, and let $d$ be the metric on $X$ induced by the Euclidean metric on $Y$.
Then $d$ cannot be extended to $\overline X$, and $(\phi|_{Y-\{0\}})^{-1}:X\to\R^n$ is a compactification.
Moreover, any semialgebraic metric on $X$ non-extendable to $\overline X$ has the same phenomenon as $d$.
Hence this example indicates that we may prove Theorem 1 by something like blows up.\par
The condition in Theorem 1 is not interesting in the non-definable case.
Let $X$ be contained and bounded in $\R^n$, and $d$ a metric on $X$.
Then $d$ is extended to a pseudo-metric to $\overline X$, i.e., $d:X\times X\to\R$ is extended to $\overline X\times\overline X\to\R$ as a continuous map if and only if $\lim_{t,t'\to0}d(\gamma(t),\gamma(t'))=0$ for any continuous curve $\gamma:(0,\,1]\to X$ such that $\gamma(t)$ converges to some point of $\overline X-X$ as $t\to0$.
However, this condition has no bearing on compactification of $X$.
Note that there is no $X$ of positive dimension such that $\lim_{t,t'\to0}d(\gamma(t),\gamma(t'))=0$ for any continuous curve $\gamma:(0,\,1]\to X$.\par
It is natural to ask whether an isometry between definable metric spaces is definable.
We give a positive answer and a negative one.
A {\it Nash set} over $R$ is the zero set of of some Nash (=\,semialgebraic $C^\infty$ smooth) function, and a {\it Nash metric} is a metric such that the graph of the function $X\times X-\Delta\ni(x,x')\to d(x,x)\in R$ is a Nash set, where $\Delta$ denotes the diagonal of $X$.
A {\it Nash metric space} is a Nash set with a Nash metric.
A typical example is an algebraic set with the Euclidean metric.\vskip1mm\noindent
{\bf Remark.} {\it An isometry between Nash metric spaces over $R$ is semialgebraic.}\vskip1mm
The reason is the following:
Let $f:X\to Y$ be an isometry between Nash metric spaces, and let $d$ denote the metrics.
We assume that $X$ is connected for simplicity of notation.
Then there are points $x_1,\ldots,x_k$ in $X$, $k=\dim X+1$, such that the map $\phi_X:X\ni x\to(d(x,x_1),\ldots,d(x,x_k))\in R^k$ is injective (but not necessarily homeomorphic), i.e., $x$ is uniquely determined by its distances to $x_1,\ldots,x_k$, and each $y\in Y$ is also so by the distances to $f(x_1),\ldots,f(x_k)$.
(This statement follows from the next note and the fact below that is proved by the identity theorem as in the proof of (2) on page 68 in \cite{F-S}.
We note that for each $x_0\in X$, the function $X\ni x\to d(x,x_0)\in R$ is nowhere locally constant by the identity theorem.
Let $\phi:X\to R^l$ be a map whose graph is a semialgebraic set.
Set $\max_{y\in R^l}\dim\phi^{-1}(y)=m$.
Then the fact is that for some point $x_0$ in $X$, the map $(\phi,\phi_0):X\ni x\to(\phi(x),d(x,x_0))\in R^{l+1}$ is injective if $m=0$ and $\dim(\phi,\phi_0)^{-1}(y)<m$ for any $y\in R^{l+1}$ otherwise.)
Hence for the proof that $f$ is semialgebraic, we can replace $X$ and $Y$ with $\Ima\phi_X$ and $\Ima\phi_Y$, respectively, because $\phi_X$ and $\phi_Y$ are semialgebraic.
By the assumption that $f$ is metric preserving, $\Ima\phi_X$ coincides with $\Ima\phi_Y$, and the map $\Ima\phi_X\to\Ima\phi_Y$ induced by $f$ is the identity map.
Therefore, $f$ is automatically semialgebraic.\par
In the same way we can prove that an isometry between real analytic metric spaces in the same sense as above is subanalytic, and it is clear that an isometry between definable sets with Euclidean metrics is definable.
However, this is not the case in general as follows:\vskip1mm\noindent
{\bf Theorem 2.} {\it There are two compact semialgebraic metric spaces over $\R$ which are isometric but not semialgebraically isometric.}\vskip1mm
The following fact shows how compact definable metric spaces are topologically simple:
Let $X$ be a compact definable metric space over $R$.
Then we can assume that $X$ is a polyhedron by the triangulation theorem of definable sets (Theorem II.2.1 in \cite{S1}).
The fact is that there exists a definable homeomorphism $\tau$ of $X\times X$ of the form $\tau(x,x')=(\tau_1(x),\tau_2(x,x'))$ such that $d\circ\tau$ is PL, i.e, piecewise linear, on $X\times X$, which we can prove as in the proof the triangulation theorem of definable continuous functions (Theorem II.3.1 in \cite{S1}).
However, this does not mean that there exists a definable homeomorphism $\tau_1$ of $X$ such that the map $X\times X\in(x,x')\to d(\tau_1(x),\tau_1(x'))\in R$ is PL, i.e., $X$ is isometric to a PL metric space.
For example, a closed interval in $R$ cannot be piecewisely linearized for $d$ defined by $d(x,x')=|x-x'|^{1/2}$.\par
For the study of topological properties of definable sets, we cannot expect to use algebra.
Hence the exact blow up is useless.
We introduce a version and use it in the proof of Theorem 1 as the key.
By the same reason we cannot obtain a resolution of singularities of definable sets by blows up.
On the other hand, we can prove it by our version as follows:\vskip1mm\noindent
{\bf Theorem 3 (Resolution of singularities of definable sets).} {\it Let $X$ be a definable set of dimension $k$.
Let $r$ be a nonnegative integer.
Then there exist a definable $C^r$ manifold $Y$ possibly with boundary and a definable $C^0$ map $\psi:Y\to X$ such that $\psi|_{\Int Y}$ is a $C^r$ diffeomorphism (homeomorphism if $r=0$) onto the set of definably $C^r$ smooth points of $X$ of dimension $k$, i.e., the subset of $X$ of points where the germ of $X$ is definably $C^r$ diffeomorphic (homeomorphic if $r=0$) to the germ of $R^k$ at 0, where $\Int Y$ denotes the interiors of $Y$.}\vskip1mm 
We thank Michel Coste.
He posed a question on Theorem 1, and we discussed it.
\section{Preparations}
Consider the case $X=(0,\,1]\times[-1,\,1]^{n-1}$, to which general $X$ is reduced in the proof of Theorem 1, and assume the condition in Theorem 1, i.e., $\lim_{t,t'\to 0}\allowbreak d(\gamma(t),\gamma(t'))=0$ for any definable continuous curve $\gamma:(0,\,1]\to X$.
We call a point $x$ of the boundary $\{0\}\times [-1,\,1]^{n-1}$ {\it bad} if there exist definable continuous curves $\gamma,\gamma':(0,\,1]\to X$ such that $\lim_{t\to0}\gamma(t)=\lim_{t\to 0}\gamma'(t)=x$ and $\lim_{t\to0}d(\gamma(t),\gamma'(t))\not=0$.
We say that $\gamma$ is of the {\it standard form} if $\gamma(t)$ is of the form $(t,\gamma_2,\ldots,\gamma_n(t))$.
(We define a {\it bad point} of $\overline X-X$ in the same way also in the case where $X$ is a general definable metric space.) 
In the proof of Theorem 1 we use the following three lemmas:\vskip1mm\noindent
{\bf Lemma 1.} {\it The set of bad points is definable and of dimension$\,<n-1$.
(If $X$ is a general definable metric space, the same statement holds, although $\overline X-X$ is not necessarily of dimension $n-1$.)}\vskip1mm\noindent
{\it Proof.} We prove only the first statement because the second is proved in the same way.
For a definable continuous curve $\gamma=(\gamma_1,\ldots,\gamma_n):(0,\,1]\to X$ with 
$\lim_{t\to 0}\gamma(t)\in\{0\}\times [-1,\,1]^{n-1}$, $\gamma_1$ is $C^1$ smooth and injective on $(0,\,t_0]$ for some small $t_0$.
We can replace $\gamma$ with $\gamma\circ\gamma^{-1}_1$ and assume that $\gamma$ is always of the standard form.
Indeed, for definable continuous curves $\gamma,\gamma':(0,\,1]\to X$ with $\lim_{t\to0}\gamma(t)\in\{0\}\times [-1,\,1]^{n-1}$ and $\lim_{t\to0}\gamma'(t)\in\{0\}\times [-1,\,1]^{n-1}$, we have $d(\gamma'\circ\gamma_1^{\prime-1}(t),\gamma\circ\gamma_1^{-1}(t))\le d(\gamma'\circ\gamma_1^{\prime-1}(t),\gamma'(t))+d(\gamma'(t),\gamma(t))+d(\gamma(t),\gamma\circ\gamma_1^{-1}(t))$.
Hence $\lim_{t\to 0}d(\gamma'\circ\gamma_1^{\prime-1}(t),\gamma\circ\gamma_1^{-1}(t))=0$ if and only if $\lim_{t\to 0}d(\gamma'(t),\gamma(t))=0$ by the hypothesis on $d$.
We write the continuous extension of $\gamma$ to $[0,\,1]$ as $\overline\gamma$.\par
First the set of bad points is definable for the following reason:
Set
\begin{align*}S=\{(t,x_2,\ldots,x_n,y_2,\ldots,y_n,u,v)\in(0,\,1]\times[-1,\,1]^{2n-2}\times R^2:\ \,\\
u=|(x_2-y_2,\ldots,x_n-y_n)|,\,v=d((t,x_2,\ldots,x_n),(t,y_2,\ldots,y_n))\}.\end{align*}
Then $S$ is definable and closed in $(0,\,1]\times[-1,\,1]^{2n-2}\times R^2$, and the set of bad points is the image of $\overline S\cap\{0\}\times[-1,\,1]^{2n-2}\times\{0\}\times(0,\,\infty)$ under the projection of $R^{2n+1}$ onto the first $n$ factors.
The reason is that there is a definable continuous curve in $\overline S$ starting from any point $a$ in $\overline S\cap\{0\}\times[-1,\,1]^{2n-2}\times\{0\}\times(0,\,\infty)$ and moving in $S$ by the curve selection lemma ((II.1.12) in \cite{S1}) and the curve induces $\gamma$ and $\gamma'$ which show that the image of $a$ under the same projection is a bad point.
Conversely, for a bad point $a_1$, we have the $\gamma$ and $\gamma'$ whose limit is $a_1$.
These $\gamma$ and $\gamma'$ define a definable continuous curve in $\overline S$ starting from a point $a$ and moving in $S$.
Then the image of $a$ under the projection is $a_1$.
Hence the image of $\overline S\cap\{0\}\times[-1,\,1]^{2n-2}\times\{0\}\times(0,\,\infty)$ is the set of bad points, and the set is definable.\par
Next we see that the set is of dimension $<n-1$ by reductio ad absurdum.
Assume the contrary.
Then by shrinking $[-1,\,1]^{n-1}$, we can suppose that $\{0\}\times(-1,\,1)^{n-1}$ (but not $\{0\}\times[-1,\,1]^{n-1}$) consists only of bad points and there are a definable continuous function $\rho$ on $[0,\,1]$ with zero set $\{0\}$ and definable continuous curves $\overline\gamma_x,\overline\gamma'_x:[0,\,1]\to X$ for each $x\in[-t_0,\,t_0]^{n-1}$ for some $t_0>0$ such that $\overline\gamma_x(0)=\overline\gamma'_x(0)=(0,x),\ |\overline\gamma_x(t)-(t,x)|\le\rho(t),\ |\overline\gamma'_x(t)-(t,x)|\le\rho(t)$ and $\lim_{t\to0}d(\overline\gamma_x(t),\overline\gamma'_x(t))\not=0$, i.e., $\rho$ is so large that for all $x\in[-t_0,\,t_0]^{n-1}$, $|\overline\gamma_x(t)-(t,x)|$ and $|\overline\gamma'_x(t)-(t,x)|$ functions in the variable $t$, are uniformly bounded by $\rho(t)$.
Here the reason of existence of such a $\rho$ is that if we remove a smaller dimensional definable subset from the set of bad points then we can choose $\gamma_x$ and $\gamma'_x$ for each $x\in[-t_0,\,t_0]^{n-1}$ so that the maps $(0,\,1]\times[-t_0,\,t_0]^{n-1}\ni(t,x)\to\gamma_x(t),\gamma'_x(t)\in X$ are definable, continuous and definably continuously extended to $[0,\,1]\times(-1,\,1)^{n-1}$.
For simplicity of notation and without loss of generality we assume $\rho(t)=t$.\par
Consider the case $n=2$.
For each $(t,x)\in(0,\,1]\times[-1,\,1]$, set $\alpha(t,x)=\sup_{0<t'\le t}d((t,x),(t',x))$ the supremum of the metrics of $(t,x)$ and points on the half-open segment with ends $(0,x)$ and $(t,x)$.
Then $\alpha$ is a definable function but not necessarily continuous, $\lim_{t\to0}\alpha(t,x)=0$ for each $x\in[-1,\,1]$ by the hypothesis on $d$, and hence shrinking $[-1,\,1]$ of $(0,\,1]\times[-1,\,1]$, we assume that $\alpha$ is continuous and, moreover, extended to a continuous function $\overline\alpha$ on $[0,\,1]\times[-1,\,1]$.
By the same reason we suppose also that the function $\beta$ on $(0,\,1/2]\times[-1,\,1/2]$, defined by $\beta(t,x)=\sup_{0<t'\le t}d((t,x),(t',x+t-t'))$ the supremum of the metrics of $(t,x)$ and points on the half-open segment with ends $(0,x+t)$ and $(t,x)$, is definable, continuous and extended to a continuous function $\overline\beta$ on $[0,\,1/2]\times[-1,\,1/2]$.
Clearly $\overline\alpha=\overline\beta=0$ on $\{0\}\times[-1,\,1/2]$.
Hence for any $\epsilon>0$ there exists small $t_0>0$ such that $\alpha\le\epsilon/3$ and $\beta\le\epsilon/3$ on $[0,\,t_0]\times[-1,\,1/2]$, which implies $d((t,x),(t',x'))\le\epsilon$ for $(t,x),\,(t',x')\in(0,\,t_0]\times[-1,\,1/2]$.
Indeed, if $x'>x$ then
\begin{align*}d((t,x),(t',x'))\le\hskip42mm\\
d((t,x),(t_0,x))+d((t_0,x),(t_0+x-x',x'))+d((t_0+x-x',x'),(t',x'))\\
\le\alpha(t_0,x)+\beta(t_0,x)+\alpha(t_0,x')\le\epsilon.\hskip25mm\end{align*}
However, by the hypothesis that $\{0\}\times(-1,\,1)$ consists only of bad points we have definable continuous curves $\gamma,\,\gamma':\allowbreak(0,\,t_0]\to(0,\,t_0]\times[-1,\,1/2]$ such that $\lim_{t\to0}\gamma(t)=\lim_{t\to0}\gamma'(t)$ and $\lim_{t\to0}d(\gamma(t),\gamma'(t))\allowbreak\not=0$, which is a contradiction.
Thus the case $n=2$ is proved.\par
Let $n>2$.
We reduce to the case essentially $n=2$.
By the assumption that points of $\{0\}\times(-1,\,1)^{n-1}$ are all bad and by the above proof that the set of bad points is definable, there exist $c>0$ and a compact definable subset $C$ of $(-1,\,1)^{n-1}$ of dimension $n-1$ such that $\lim_{t\to 0}d(\gamma_x(t),\gamma_x'(t))\ge c$ for each $x\in C$ and for some definable continuous curves $\gamma_x,\,\gamma_x':(0,\,1]\to X$ with $\lim_{t\to 0}\gamma_x(t)=\lim_{t\to 0}\gamma_x'(t)=(0,x)$.
Here we suppose the following three conditions:
First $C=[-1/2,\,1/2]^{n-1}$ always for simplicity of notation, even if we shrink $C$ in the argument below.
Secondly, there exist definable continuous maps $\gamma,\,\gamma':(0,\,1]\times C\to X$ such that $\gamma_x(t)=\gamma(t,x)$ and $\gamma'_x(t)=\gamma'(t,x)$.
The $\gamma_x(t)$ and $\gamma'_x(t)$ are of the standard form.
Thirdly, $\gamma_x(t)=(t,x)$ and $\gamma'_x(t)=(t,x_2,\ldots,x_{i-1},\gamma'_{x,i}(t),x_{i+1},\ldots\allowbreak,x_n)$ for some $i$ for the following reason:\par
Set
$$\gamma^i_x(t)=(t,\gamma_{x,2}(t),\ldots,\gamma_{x,i}(t),\gamma'_{x,i+1}(t),\ldots,\gamma'_{x,n}(t))\ \text{for }i=1,\ldots,n.$$
Then $\gamma^1_x=\gamma'_x,\ \gamma^n_x=\gamma_x$ and $\lim_{t\to0}\gamma^i_x(t)=(0,x)$.
Consider the pairs $(\gamma^{i-1}_x,\gamma^i_x),\ i=2,\ldots,n$.
At least one of them, say, $(\gamma^{i-1}_x,\gamma^i_x)$ has the same property as $(\gamma_x,\gamma'_x)$, i.e., $\lim _{t\to0}d(\gamma^{i-1}_x(t),\gamma^i_x(t))\ge c/(n-1)$.
Hence we can replace $(\gamma_x,\gamma'_x)$ by $(\gamma^{i-1}_x,\gamma^i_x)$, and assume that $\gamma_x(t)$ and $\gamma'_x(t)$ are of the forms, respectively,
\begin{align*}(t,\gamma_{x,2}(t),\ldots,\gamma_{x,i-1}(t),\gamma_{x,i}(t),\gamma_{x,i+1}(t),\ldots,\gamma_{x,n}(t))\ \text{and}\\
(t,\gamma_{x,2}(t),\ldots,\gamma_{x,i-1}(t),\gamma'_{x,i}(t),\gamma_{x,i+1}(t),\ldots,\gamma_{x,n}(t)).\quad\ \end{align*}\par
We need to further simplify the forms.
For this we need the property that the map $(0,\,1]\times C\ni(t,x)\to\gamma_x(t)\in(0,\,1]\times(-1,\,1)^{n-1}$ is an embedding.
Let $D_{t_0}$ denote the set of points of $(0,\,t_0]\times C$ where $\gamma$ is not $C^1$ regular for $0<t_0<1$.
Assume that $D_{t_0}$ is of dimension $n$ for any small $t_0$.
Then $\gamma(D_{t_0})$ is of dimension $<n$, and by shrinking $C$ and decreasing $t_0$, we can assume that $D_{t_0}=(0,\,t_0]\times C$.
But $\overline{\gamma(D_{t_0})}-\gamma(D_{t_0})$ contains $\{0\}\times C$ since every set and map are definable, which contradicts the fact $\dim\gamma(D_{t_0})<n$.
Hence $D_{t_0}$ is of dimension $<n$ for some small $t_0$, and we can assume $D_{t_0}=\emptyset$.
Moreover, by similar arguments we see the local covering number of the immersion $\gamma$ at each point of the image near $\{0\}\times(-1/2,\,1/2)^{n-1}$ is one, and hence the continuous extension of $\gamma|_{(0,\,t_0]\times C}$ to $[0,\,t_0]\times C$ is a continuous embedding into $[0,\,1]\times(-1,\,1)^{n-1}$.
Hence $\gamma^{-1}\circ\gamma'_x:(0,\,t_0]\to X$ is a well-defined continuous curve if we shrink $C$ and decrease $t_0$, and we can replace $\gamma_x$ and $\gamma'_x$ by the maps $(0,\,t_0]\ni t\to(t,x)\in X$ and $\gamma^{-1}\circ\gamma'_x$ respectively.
Then $\gamma_x(t)$ and $\gamma'_x(t)$ are of the required forms, i.e., the third condition holds.\par
We can regard $\gamma$ and $\gamma'$ satisfying the third condition as curves in $(0,\,t_0)\times R$.
Therefore, Lemma 1 is proved also in the case $n>2$.\qed\vskip2mm
Let $X$ and $d$ be as above.
We call a bad point $x$ {\it r-bad} for a positive definable continuous function $r=r(t)$ on $(0,\,1]$ if there exist definable continuous curves $\gamma,\,\gamma':(0,\,1]\to X$ of the standard form such that $\lim_{t\to 0}\gamma(t)=\lim_{t\to 0}\gamma'(t)=x$, $\lim_{t\to 0}d(\gamma(t),\gamma'(t))\not=0$ and $\{|\gamma(t)-\gamma'(t)|/r(t):t\in(0,\,1]\}$ is bounded.
(If $X$ is a general definable metric space, a curve of the standard form is meaningless.
In this case we replace the condition that $\{|\gamma(t)-\gamma'(t)|/r(t)\}$ is bounded by one that $\{|\gamma(s)-\gamma'(s')|/r(t):s,s',t\in(0,\,1],\,|\gamma(s)|=t,\,|\gamma'(s')|=t\}$ is bounded or, equivalently, that $|\gamma\circ|\gamma|^{-1}(t)-\gamma'\circ|\gamma'|^{-1}(t)|/r(t)$ is bounded, where $|\gamma|$ denotes the function $(0,\,1]\ni t\to|\gamma(t)|\in R$, and we call the point $\lim_{t\to0}\gamma(t)$ {\it$r$-bad}.)
We say that $x$ {\it is not an $r$-bad point} if it is not a bad point or it is a bad point but not an $r$-bad point.
Then the second lemma is the following:\vskip1mm\noindent
{\bf Lemma 2.} {\it There are no $r$-bad points for some $r$.
This is weaker than (actually equivalent to) the condition $|\gamma(t)-\gamma'(t)|>r(t),\ t\in(0,\,t_0]$, for some $r$, for any two definable continuous curves $\gamma$ and $\gamma'$ of the standard form such that $\lim_{t\to0}\gamma(t)=\lim_{t\to0}\gamma'(t)\in\{0\}\times[-1,\,1]^{n-1}$ and $\lim_{t\to0}d(\gamma(t),\gamma'(t))\not=0$ and for some $t_0\in(0,\,1]$.
(The same statement holds for a general definable metric space if the condition $|\gamma(t)-\gamma'(t)|>r(t)$ is replaced by one $|\gamma\circ|\gamma|^{-1}(t)-\gamma'\circ|\gamma'|^{-1}(t)|>r(t)$.)}\vskip1mm\noindent
{\it Proof.} We treat only the case $X=(0,\,\infty)\times[-1,\,1]^{n-1}$.
Let $r=r(t)$ be a positive definable continuous function on $(0,\,1]$.
We first prove that the set of $r$-bad points is definable.
Let $S$ be the set defined in the proof of Lemma 1, and set
\begin{align*}S_r=\{(t,x_2,\ldots,x_n,y_2,\ldots,y_n,u,v,w)\in S\times R:u\le w r(t)\},\\
S'_r=\cup_{w\in R}(S_r\cap R^{2n+1}\times[0,\,w]).\hskip20mm\end{align*}
Then the set of $r$-bad points is the image of $\overline{S'_r}\cap\{0\}\times[-1,\,1]^{2n-2}\times\{0\}\times(0,\,\infty)\times R$ under the projection of $R^{2n+2}$ onto the first $n$ factors by the same reason as that the set of bad points is the image of $\overline S\cap\{0\}\times[-1,\,1]^{2n-2}\times\{0\}\times(0,\,\infty)$.
Hence the set of $r$-bad points is definable.\par
If we fix $\gamma$ and $\gamma'$, we of course have $r$ such that $|\gamma(t)-\gamma'(t)|>r(t)$.
Next we will choose $r$ independently of $\gamma$ an $\gamma'$ by a method of quantifier elimination.
For each $v\in(0,\,\infty)$, let $T(v)$ denote the subset of $(0,\,1]\times[0,\,1]$ consisting of points $(t,u)$ such that $t\le v/2$ and $(t,x_2,\ldots,x_n,y_2,\ldots,y_n,u,v)\allowbreak\in S$ for some $(x_2,\ldots,x_n,y_2,\ldots,y_n)\in[-1,\,1]^{2n-2}$.
Set $T=\cup_{v\in(0,\,\infty)}T(v)\times\{v\}$.
Then for each $v\in(0,\,\infty)$, $\overline{T(v)}\cap[0,\,1]\times\{0\}\subset\{(0,0)\}$ because $T(v)$ is contained and closed in $(0,\,1]\times[0,\,1]$, because $T(v)\cap(0,\,1]\times\{0\}=\emptyset$ and hence because some definable neighborhood of $(0,\,1]\times\{0\}$ in $[0,\,1]\times[0,\,1]$ does not touch $T(v)$.
Hence $\overline T\cap[0,\,\infty)\times\{0\}\times[0,\,\infty)\subset\{0\}\times\{0\}\times[0,\,\infty)$ because $t\le v/2$ for $(t,u,v)\in\overline T$.
Therefore, $q(\overline T)$, the image under the projection $q$ of $[0,\,1]^2\times[0,\,\infty)$ forgetting the last factor, intersects with $[0,\,1]\times\{0\}$ at most at the origin.
Clearly $q(\overline T)$ is definable and closed in $R^n$.
Consequently, there exists a non-negative definable continuous function $r$ on $[0,\,1]$ such that $r^{-1}(0)=\{0\}$ and $q(\overline T)\cap\{(t,u)\in[0,\,1]^2:r(t)\ge u\}\subset\{(0,0)\}$, i.e., $q(T)\subset\{(t,u)\in(0,\,1]\times[0,\,1]:r(t)<u\}$.
The last inclusion means that $|\gamma(t)-\gamma'(t)|\ge r(t),\ t\in(0,\,t_0]$, for any definable continuous curves $\gamma,\gamma':(0,\,1]\to(0,\,1]\times[-1,\,1]^{n-1}$ of the standard form such that $\lim_{t\to0}\gamma(t)=\lim_{t\to0}\gamma'(t)$ and $\lim_{t\to0}d(\gamma(t),\gamma'(t))\not=0$ if exist and for some $t_0\in(0,\,1]$.
Indeed, for such $\gamma$, $\gamma'$ and small $t>0$ we have
$$(t,|\gamma(t)-\gamma'(t)|,d(\gamma(t),\gamma'(t)))\in T,\text{ hence }(t,|\gamma(t)-\gamma'(t)|)\in q(T)\text{ and }r(t)<|\gamma(t)-\gamma'(t)|.$$
Hence there is no $r'$-bad (not $r$-bad) point, where $r'$ is the function defined by $r'(t)=r(t)/t$.\qed\vskip2mm\noindent
{\bf Remark.} {\it For each $c\in[-1,\,1]^{n-k}$, set $X_c=(0,\,1]\times[-1,\,1]^{k-1}\times\{c\}$ and $d_c=d|_{X_c}$, and define $c$-bad points of $\{0\}\times[-1,\,1]^{k-1}\times\{c\}$ and $r$-$c$-bad points for the above $r$ by $d_c$ likewise bad points and $r$-bad points by $d$.
Then a $c$-bad point is a bad point, and an $r$-$c$-bad point is an $r$-bad point.}\vskip2mm
If a manifold has boundary then we call it a {\it manifold with boundary}.
A $C^1$ {\it manifold with corners} is a $C^1$ manifold in which each point has a neighborhood $C^1$ diffeomorphic to $R\times\cdots\times R\times[0,\,\infty)\times\cdots\times[0,\,\infty)$.
We may call a $C^1$ manifold with boundary a $C^1$ manifold with corners.
By the symbols $\partial$ and $\Int$ we mean the boundary and the interior of a $C^1$ manifold with corners but not those as a subset of its ambient Euclidean space.
For definable sets $X\subset Y$, let $\bdry X$ denote the boundary of $X$ in $Y$ as a topological subspace when $Y$ is clear by the context.
Let $Y$ be a $C^1$ manifold possibly with corners, and $Z_1$ and $Z_2$ $C^1$ submanifolds possibly with corners.
We say that $Z_1$ {\it intersects transversally} $Z_2$ in $Y$ if the tangent spaces of $Z_1$ and $Z_2$ at any point of $Z_1\cap Z_2$ span the tangent space of $Y$ at the point as linear spaces.
For example, the segment $[0,\,1]\times\{0\}$ intersects transversally the one $\{0\}\times[0,\,1]$ in $[0,\,1]^2$.\par
Let $X$ be a $C^1$ manifold with corners, and $\{X_0,\ldots,\allowbreak X_{\dim X}\}$ be the stratification of $X$ into definable $C^1$ manifolds such that $X_i$ is the set of definably $C^1$ smooth points of $X-\cup_{i'=0}^{i-1}X_{i'}$ of dimension $i,\ i=0,\ldots,\dim X$.
Note $X_{\dim X}=\Int X$.
We call $\{X_0,\ldots,\allowbreak X_{\dim X}\}$ the {\it canonical stratification} of $X$.
We make\vskip1mm\noindent
{\it the face assumption that the closure of each connected component of each $X_i$ is a definable $C^1$ manifold possibly with corners.}\par
We call the closure a {\it face} of $X$, a face a {\it proper face} if it is not a connected component of $X$.
Note that a face of a face is a face.
We have in mind an example of such a manifold with corners: $[0,\,\infty)\times R^{n-1}\cap N,\ n>2$, where $N=\{x\in R^n:|x|\ge1\}$, whose proper faces are $[0,\,\infty)\times R^{n-1}\cap\partial N$, $\{0\}\times R^{n-1}\cap N$ and $\{0\}\times R^{n-1}\cap\partial N$.
An example of a $C^1$ manifold $X$ with corners, of dimension 2 and with the canonical stratification $\{X_0,X_1,X_2\}$ such that $X_1$ is connected and $\overline{X_1}$ is not a $C^1$ manifold possibly with corners is $\{(x,y)\in\R^2:x y\ge(x^2+y^2)^2,\,x\ge0\}$.
Recall that given a $C^1$ manifold $\tilde Y$, its $C^1$ submanifold $Y$ with corners satisfying the face assumption and another $C^1$ submanifold $Z$, if $Z$ intersects transversally each face of $Y$ in $\tilde Y$ then $Y\cap Z$ is a $C^1$ manifold possibly with corners satisfying the face assumption and $\partial(Y\cap Z)=\partial Y\cap Z$.\vskip1mm\noindent
{\bf Definition 1.} The idea of the proof of Theorem 1 is to use blows up.
However, we exchange blows up with hole-blows-up, which are convenient for our problem, because for the blow up $\psi:Y\to R^n$ of $R^n$ with center 0, we cannot characterize $\psi^{-1}(0)$ in $Y$, while if $\psi:Y\to R^n$ is the hole-blow-up of $R^n$ with center 0, then $Y$ is a definable $C^1$ manifold with boundary and $\psi^{-1}(0)$ is the boundary.
A {\it hole-blow-up} of $R^n$ with center 0 is the pair of a definable $C^1$ manifold $N$ with boundary and a definable $C^1$ map $\psi:N\to R^n$ such that $N=\{x\in R^n:|x|\ge\epsilon\}$ for some small $\epsilon>0$, $\psi|_{N-\partial N}$ is a diffeomorphism onto $R^n-\{0\}$ and $\psi(x_1,\ldots,x_n)$ is of the form $(x_1,x_1x_2,\ldots,x_1x_n)$ for some local coordinate systems $(x_1,\ldots,x_n)$ of $N$ at each point of $\partial N$ and of $R^n$ at 0.
Such a $\psi$ is defined, e.g., by $\psi(x)=(1-\epsilon/|x|)x$.
Note that the blow up of $R^n$ with center 0 is the quotient space of $N$ where two points $x$ and $-x$ in $\partial N$ are identified.
For a definable $C^1$ manifold $X$ of dimension $k$ contained in $R^n$ and containing 0, we define a {\it hole-blow-up} of $X$ with center 0 to be the map $\psi:X\cap N\to X$ for small $\epsilon>0$ through a definable $C^1$ diffeomorphism $\phi:\{x\in X:|x|\le2\epsilon\}\to\{x\in R^k:|x|\le2\epsilon\}$ with $|\phi(x)|=|x|$.
Here we choose $\psi$ to be the identity map on $\{x\in X:|x|\ge3\epsilon\}$ and so that its restriction to $X\cap\Int N$ is a diffeomorphism onto $X-\{0\}$.
If $X$ is a definable $C^1$ manifold contained in $R^n$, $\overline X$ is a definable $C^1$ manifold with boundary $\overline X-X$ and $\overline X-X$ contains 0, then we define a {\it hole-blow-up} $\psi:\overline X\cap N\to\overline X$ of $X$ with center 0 so that $\psi(X\cap N)=X\cup\{0\}$ by using a local coordinate system of $\overline X$ at 0 in the same way.
Note that $\overline X\cap N$ is a definable $C^1$ manifold with corners.
Other notes are the following:\vskip1mm\noindent
{\it Consider the case of $X=(0,\,\infty)\times R^{n-1}$ with a definable metric $d$ on $X$.
Let $\psi:\overline X\cap N\to\overline X$ be a hole-blow-up of $X$ with center 0.
Assume that $0$ is not a $t$-bad point, i.e., for any definable continuous curves $\gamma,\gamma':(0,\,1]\to X-\{0\}\times R^{n-1}$ of the standard form such that $\lim_{t\to0}\gamma(t)=\lim_{t\to0}\gamma'(t)=0$ and $|\gamma(t)-\gamma'(t)|/t$ is bounded, the $d(\gamma(t),\gamma'(t))$ converges to 0 as $t\to0$.
(We do not assume that $d$ is continuously extended to $\overline X-\{0\}$.)
Then the metric $d_\psi$ on $X\cap\Int N$ defined by $d_\psi(x,x')=d(\psi(x),\psi(x'))$ is continuously extended to $(X\cap\Int N)\cup(\overline X\cap\partial N)$.\par
Let $X=R^n-\{0\}\times R^k$ with a definable metric $d,\ 0\le k<n-1$.
(If $k=n-1$, this note is similar to the above one but false.)
Let a function $r$ on $R^n$ be defined by $r(t)=|x|$.
Let $\psi:N\to R^n$ be a hole-blow-up of $R^n$ with center 0 such that $\psi(X\cap N)=X$ as above.
Assume that 0 is not an $r$-bad point.
Then the metric $d_\psi$ on $X\cap\Int N$ defined by $d$ and $\psi$ as above is continuously extended to $(X\cap\Int N)\cup(\overline X\cap\partial N)$.}\par
We say that $d_\psi$ is {\it defined by} $d$ and $\psi$.
The continuous extension in the notes means that the function $(X\cap\Int N)^2\ni(x,x')\to d_\psi(x,x')\in R$ is continuously extended to $((X\cap\Int N)\cup(\overline X\cap\partial N))^2$ or, equivalently, $d_\psi$ is extended to a pseudo metric on $(X\cap\Int N)\cup(\overline X\cap\partial N)$ but not to a metric there.\par
{\it Proof of the notes.} We prove the first note by reductio ad absurdum.
Assume the contrary.
Then there exist definable continuous curves $\gamma,\gamma':(0,\,1]\to X\cap\Int N$ such that $\lim_{t\to0}\gamma(t)=\lim_{t\to0}\gamma'(t)\in\overline X\cap\partial N$ and $\lim_{t\to0}d_\psi(\gamma(t),\gamma'(t))\not=0$.
There are two cases: $\lim_{t\to0}\gamma(t)\in X\cap\partial N$ or $\lim_{t\to0}\gamma(t)\in\overline X\cap\partial N-X$.
Consider the first case.
Choose local definable $C^1$ coordinate systems $(t,x_2,\ldots,x_n)$ of $N$ at the point $\lim_{t\to0}\gamma(t)$ and $(y_1,\ldots,y_n)$ of $R^n$ at the origin so that $N=\{(t,x_2,\ldots,x_n):t\ge0\}$, $\lim_{t\to0}\gamma(t)=(0,\ldots,0)$, the origin in $R^n$ is $(0,\ldots,0)$ of the coordinate system and $\psi(t,x_2,\ldots,x_n)=(t,t x_2,\ldots,t x_n)$.
Then we can suppose that $\gamma$ and $\gamma'$ are of the standard form, i.e., $\gamma(t)$ and $\gamma'(t)$ are of the form $(t,\gamma_2(t))$ and $(t,\gamma'_2(t))$, respectively, as in the proof of Lemma 1 since $d_\psi$ on $X\cap\Int N$ satisfies the condition in Theorem 1.
It follows $\psi\circ\gamma(t)=(t,t\gamma_2(t))$ and $\psi\circ\gamma'(t)=(t,t\gamma'_2(t))$, and we have
$$d_\psi(\gamma(t),\gamma'(t))=d(\psi\circ\gamma(t),\psi\circ\gamma'(t))=d\big((t,t\gamma_2(t)),(t,t\gamma'_2(t))\big).$$
\vskip-2mm\noindent Hence\vskip-5mm
$$\lim_{t\to0}(t,t\gamma_2(t))=\lim_{t\to0}(t,t\gamma'_2(t))=0,\ \ \lim_{t\to0}d\big((t,t\gamma_2(t)),(t,t\gamma'_2(t))\big)\not=0$$
and $|(t,t\gamma_2(t))-(t,t\gamma'_2(t))|/t\,(=|\gamma_2(t)-\gamma'_2(t)|)$ is bounded, which means that 0 in $R^n$ is a $t$-bad point with respect to $d$.
This contradicts the hypothesis that 0 is not a $t$-bad point.
In the case $\lim_{t\to0}\gamma(t)\in\overline X\cap\partial N-X$, we modify $\gamma$ and $\gamma'$ so that $\gamma(t)=(t,\gamma_2(t))$ and $\gamma'(t)=(t,\gamma'_2(t))$ for some local definable $C^1$ coordinate systems $(t,x_2,\ldots,x_n)$ of $N$ at the point $\lim_{t\to0}\gamma(t)$ and $(y_1,\ldots,y_n)$ of $R^n$ at 0 such that $\lim_{t\to0}\gamma(t)=(0,\ldots,0)$, $X\cap N=\{(t,x_2,\ldots,x_n)\in R^n:t\ge0,\,x_2>0\},\ X=\{(y_1,\ldots,y_n)\in R^n:y_2>0\}$ and $\psi(t,x_2,\ldots,x_n)=(t,t x_2,\ldots,t x_n)$.
Then $\psi\circ\gamma(t)$ and $\psi\circ\gamma'(t)$ are of the form $(t,t\gamma_2(t))$ and $(t,t\gamma'_2(t))$ respectively.
Any of them is not of the standard form since $X=\{y_2>0\}$.
However, the same contradiction is obtained because the set $\{|\gamma\circ|\gamma|^{-1}(t)-\gamma'\circ|\gamma'|^{-1}(t)|/t:t\in[0,\,1]\}$ is bounded.
The second note follows in the same way.\qed\vskip1mm
We continue the definition of a hole-blow-up.
Let $X=(0,\,\infty)\times R^{n-1}$, and $\{a_1,\ldots,a_k\}$ be a finite subset of $\{0\}\times R^{n-1}$.
We define a {\it hole-blow-up} $\psi:\overline X\cap N\to\overline X$ of $X$ with center $\{a_1,\ldots,a_k\}$ by induction on $k$ as follows:
If $k=1$, it is already defined.
Assume that a hole-blow-up $\psi_{k-1}:\overline X\cap N_{k-1}\to\overline X$ of $X$ with center $\{a_1,\ldots,a_{k-1}\}$ is given so that $N_{k-1}=\{x\in R^n:|x-a_l|\ge\epsilon_l,\,l=1,\ldots,k-1\}$ and $a_k\in\Int N_{k-1}$.
Choose $\epsilon_k>0$ so small that the set $\{x\in R^n:|x-a_k|\le3\epsilon_k\}$ is contained in $N_{k-1}$.
Set $N_k=\{x\in N_{k-1}:|x-a_k|\ge\epsilon_k\}$, and naturally define a hole-blow-up $\phi_k:\overline X\cap N_k\to\overline X\cap N_{k-1}$ of $X\cap N_k$ with center $a_k$.
Then $\psi=\psi_{k-1}\circ\phi_k:\overline X\cap N_k\to\overline X$ is a hole-blow-up of $X$ with center $\{a_1,\ldots,a_k\}$.
As in the case $k=1$, if any of $a_1,\ldots,a_k$ is not a $t$-bad point, then the metric $d_\psi$ on $X\cap N$ defined by $d$ and $\psi$ is continuously extended to $(X\cap N)\cup(\overline X\cap\partial N)$.\par
Next we consider the center of positive dimension.
First a {\it hole-blow-up} of $R^n$ with center $C=\{0\}\times R^{k-1}\,(\subset R^n)$ is the pair of the definable $C^1$ manifold $N=\{(x,y)\in R^{n-k+1}\times R^{k-1}:|x|\ge\epsilon\}$ for some small $\epsilon>0$ and a definable $C^1$ map $\psi:N\to R^n$ such that $\psi|_{N-\partial N}$ is a diffeomorphism onto $R^n-C$ and $\psi(x_1,\ldots,x_n)$ is of the form $(x_1,x_1x_2,\ldots,x_1x_{n-k},x_{n-k+1},\ldots,x_n)$ for local coordinate systems $(x_1,\ldots,x_n)$ of $N$ at a point of $\partial N$ and of $R^n$ at the image of the point under $\psi$.
For $X=(0,\,\infty)\times R^{n-1}$, we define a {\it hole-blow-up} of $X$ with center $C$ to be $\psi|_{\overline X\cap N}:\overline X\cap N\to\overline X$ for a hole-blow-up $\psi$ of $R^n$ with center $C$ such that $\psi^{-1}(X)=X\cap N$.\par
We define a hole-blow-up in a general setting.
Let $X$ be a definable $C^1$ manifold possibly with corners contained in $R^n$ such that $\overline X$ is also a definable $C^1$ manifold possibly with corners satisfying the face assumption and $\overline X-X$ is a face of $\overline X$.
Let $C$ be a definable contractible $C^1$ manifold possibly with corners contained and closed in $\overline X-X$ such that $\partial C\subset\partial(\overline X-X)$, the face assumption holds and $C$ intersects transversally to each proper face of $\overline X-X$ in $\overline X-X$.
Then we have a definable $C^1$ tubular neighborhood $V$ of $C$ in $\overline X$, a definable $C^1$ submersive retraction $\pi:V\to C$ and a definable $C^1$ diffeomorphism $\rho=(\pi,\rho_2):V\to C\times\pi^{-1}(c)$ for some $c\in C$ such that the map $V-C\ni x\to(\pi(x),|x-\pi(x)|)\in C\times(0,\,\infty)$ is a proper $C^1$ submersion onto $\{(x,t)\in C\times(0,\,\infty):t<2\epsilon\circ\pi(x)\}$ for some small positive definable $C^1$ function $\epsilon$ on $C$, $\rho(V\cap X)=C\times(\pi^{-1}(c)\cap X)$, $|\rho_2(x)-c|/\epsilon(c)=|x-\pi(x)|/\epsilon\circ\pi(x)$ and $\rho_2=\id$ on $\pi^{-1}(c)$.
Hence $V=\{x\in\overline X:|x-\pi(x)|<2\epsilon\circ\pi(x)\}$.
(Existence of $\pi$ is easily shown by a definable $C^1$ partition of unity and an elementary argument.
We omit the details of the proof (see page 162 in \cite{S1} for the construction of a definable $C^1$ partition of unity).
Then existence of $\rho$ follows from the assumption that $C$ is contractible and the o-minimal version of Thom's first isotopy lemma, Theorem II.6.1 in \cite{S1}.)
Set $N_C=\overline X-\{x\in V:|x-\pi(x)|<\epsilon\circ\pi(x)\}$, which is a definable $C^1$ manifold possibly with corners satisfying the face assumption.
Let $\bdry N_C$ always be the boundary of the topological subspace $N_C$ of $\overline X$.
Note that $\bdry N_C$ is equal to $(\partial N_C-\partial \overline X\,\overline)$ and a definable $C^1$ manifold possibly with corners satisfying the face assumption.
Then the $\pi$ restricted to $\bdry N_C$ is extended to a definable $C^1$ map $\psi:N_C\to\overline X$ so that $\psi=\id$ on $\overline X-V$, $\pi\circ\psi=\pi$ on $V\cap N_C$, $\psi|_{N_C-\bdry N_C}$ is a diffeomorphism onto $\overline X-C$, $\psi(X\cap N_C)=X$ and $\psi$ is described by a local coordinate system at each point of $\bdry N_C$ as in the case $X=(0,\,\infty)\times R^{n-1}$ and $C=\{0\}\times R^{k-1}$ by existence of $\rho$.
We call $\psi:N_C\to\overline X$ a {\it hole-blow-up} of $X$ with center $C$ and such ($X,C)$ a {\it pair of definable $C^1$ manifolds possibly with corners admitting a hole-blow-up}.
An example of such a pair is $((0,\,\infty)\times[-1,\,1]^{n-1},\{0\}\times[-1,\,1]^{k-1}),\ k<n$.
(We can obviously define a hole-blow-up of $X$ with center $C$ without the assumption that $C$ is contractible or existence of $\rho$.
We assume this for simplicity of presentation.
Hence we do not assume it in the argument below.)\par
In the argument below we will shrink only $\epsilon$.
To be precise, we replace $\epsilon$ by smaller $\epsilon'$, we define $V_{\epsilon'}$ by $\epsilon'$ in the same way, we set $\pi_{\epsilon'}=\pi|_{V_{\epsilon'}}$ and we obtain a hole-blow-up $\psi_{\epsilon'}:N_{C_{\epsilon'}}\to\overline X$ by $V_{\epsilon'}$ and $\pi_{\epsilon'}$.
We say that $\psi:N_C\to\overline X$ is {\it normally modified} (to $\psi_{\epsilon'}:N_{C_{\epsilon'}}\to\overline X$).
Note the following fact:
Let $d$ be a definable metric on $X$, and $d_\psi$ the metric on $X\cap N_C-\bdry N_C$ defined by $d$ and $\psi$.
Assume the restriction of $d_\psi$ to $\pi^{-1}(x)\cap X\cap N_C-\bdry N_C$ is continuously extended to $\pi^{-1}(x)\cap((X\cap N_C)\cup\bdry N_C)$ for each $x\in C$.
Then this property depends only on $\pi$ and hence continues to hold after any normal modification of $\psi$.\vskip1mm\noindent
{\bf Definition 2.} Let $r$ be a positive integer.
Let $C$ be a definable $C^r$ manifold in $R^n$.
A {\it definable tube} at $C$ is a triple $U=(|U|,\pi,\xi)$, where $|U|$ is a definable open neighborhood of $C$ in $R^n$, $\pi:|U|\to C$ is a definable submersive $C^r$ retraction and $\xi$ is a definable nonnegative continuous function on $|U|$ such that $\xi^{-1}(0)=C$ and $(\pi,\xi)|_{|U|-C}:|U|-C\to C\times R$ is a $C^r$ submersion.
Let $\{C_i:i=1,\ldots,k\}$ be a stratification of a definable set in $R^n$ into definable $C^r$ manifolds such that $\dim C_i\le\dim C_{i'}$ if $i\le i'$.
We call $\{C_i\}$ a {\it Whitney definable $C^r$ stratification} if the following Whitney condition holds:\par
Let $a_l$ and $b_l,\ l\in\N$, be sequences in $C_{i'}$ and $C_i$ with $i<i'$, respectively, both converging to a point $b$ in $C_i$.
Assume that the sequence of the tangent spaces $T_{a_l}C_{i'},\ l\in\N$, converges to a subspace $T$ of $R^n$ in the Grassmannian $G_{n,n'}$ and the sequence of the lines in $R^n$ passing through 0 and $a_l-b_l,\ l\in\N$, converges to a line $L$ in $R^n$ in $G_{n,1}$, where $n'=\dim C_{i'}$.
Then the Whitney condition says that $L\subset T$.\par
A {\it definable tube system} for $\{C_i\}$ consists of one definable tube $U_i=(|U_i|,\pi_i,\xi_i)$ at each $C_i$.
We call a definable tube system $\{U_i\}$ for $\{C_i\}$ {\it controlled} if each pair $i$ and $i'$ with $(\overline{C_{i'}}-C_{i'})\cap C_i\not=\emptyset$, the following condition holds:
$$\pi_i\circ\pi_{i'}=\pi_i\ \ \text{and }\ \xi_i\circ\pi_{i'}=\xi_i\ \ \text{on }|U_i|\cap|U_{i'}|.$$
We easily see that $\xi_i|_{U_i\cap C_{i'}}$ is $C^r$ regular by the Whitney condition.
From this it follows that for each $i$, the sets $R^n-\cup_{i'<i}\{x\in U_{i'}:\xi_{i'}(x)<\epsilon_{i'}\}$ and $C_i-\cup_{i'<i}\{x\in U_{i'}:\xi_{i'}(x)<\epsilon_{i'}\}$ are definable $C^r$ manifolds possibly with corners for sufficiently small positive numbers $\epsilon_1,\ldots,\epsilon_k$ such that $\epsilon_1\gg\cdots\gg\epsilon_k$.\par
We know the following:
First, any definable set admits a Whitney definable $C^r$ stratification ((II.1.14) in \cite{S1}).
Secondly, there exists a controlled definable tube system for any Whitney definable $C^r$ stratification (Lemma II.6.10, ibid.).
Thirdly, we can choose a Whitney definable $C^r$ stratification of a definable set $D$ so that the subset of $D$ of definably $C^r$ smooth points of dimension $=\dim D$ is one of the strata ((II.1.10), ibid.).\vskip1mm\noindent
{\bf Definition 3.} Let $\{C_i:i=1,2,\ldots\}$ be a Whitney definable $C^r$ stratification of a definable set $X$ locally closed in its ambient Euclidean space such that $\overline{C_i}\supset C_1$ for any $i$.
Let $U=(|U|,\pi,\xi)$ be a definable tube at $C_1$, and $\epsilon$ be a sufficiently small positive number.
Let $C_0$ be a compact definable $C^r$ manifold possibly with corners contained in $\Int C_1$ and of the same dimension.
Set
$$X_0=\{x\in X\cap|U|:\pi(x)\in C_0,\xi(x)\le2\epsilon\}\ \ \text{and}\ \ X_{00}=\{x\in X_0:\xi(x)\ge\epsilon\}.$$
Then by the Whitney condition and the o-minimal version of Thom's first isotopy lemma, there is a {\it hole-blow-up} $\psi:X_{00}\to X_0$ of $X_0-C_0$ with center $C_0$, i.e., $\psi$ is a definable $C^0$ map, $\psi=\id$ on $\{x\in X_{00}:\xi(x)=2\epsilon\}$, $\pi\circ\psi=\pi$, $\xi\circ\psi=2(\xi-\epsilon)$ and $\psi$ restricted to $C_i\cap\{x\in X_{00}:\xi(x)>\epsilon\}$ is a $C^r$ diffeomorphism onto $C_i\cap X_0$ for $i>1$.
However, $\psi$ is not necessarily of the form of a blow up locally at a point $x$ of $X_{00}$ with $\xi(x)=\epsilon$.
We can extend $\psi$ to $\psi:\{x\in X\cap\pi^{-1}(C_0):\xi(x)\ge\epsilon\}\to X\cap\pi^{-1}(C_0)$ by setting $\psi=\id$ outside $X_{00}$.
Note that $\psi$ depends on $\{C_i\}$ and $U$.\par
We can define a hole-blow-up of a definable set with center any definable closed subsets (see Figure 1).
We will explain and use it in the proofs of Theorems 1 and 3.
\vskip2mm
Consider a pair $(X,C)$ of definable $C^1$ manifolds possibly with corners admitting a hole-blow-up $\psi:N_C\to\overline X$.
Let $d$ be a definable metric on $X$.
Then it is natural to ask whether the metric $d_\psi$ on $X\cap N-\bdry N$ defined by $d$ and $\psi$ is continuously extended to $(X\cap N)\cup\bdry N$.
If this is the case, Theorem 1 is easy to prove.
However, this is false except for $\dim X\le2$.
The following third lemma solve the problem:\vskip1mm\noindent
{\bf Lemma 3.} {\it Let $(X,C)$ be a pair of definable $C^1$ manifolds possibly with corners admitting a hole-blow-up (Definition 1), and let $\psi:N_C\to\overline X$ be a hole-blow-up of $X$ with center $C$.
Let $d$ be a definable metric on $X$.
Assume that each point of $C$ is not a $t$-bad point.
(We do not assume that $d$ is continuously extended to $\overline X-C$.)
Then there exists a definable closed subset $C'$ of $C$ of dimension $<\dim C$ such that the metric $d_\psi$ on $X\cap N_C-\bdry N_C$ defined by $d$ and $\psi$ is continuously extended to $(X\cap N_C)\cup\bdry N_C-\psi^{-1}(C')$.\par
Let $N_{C'}$ be the complement in $\overline X$ of a definable open neighborhood $U_{C'}$ of $C'$ in $\overline X$ such that $N_{C'}$ is a definable $C^1$ manifold with corners satisfying the face assumption and each proper face of $N_{C'}$ intersects transversally $C$ in $\overline X$.
Let $\psi_{C'}:N_{C'}\to\overline X$ be a definable continuous map such that its restriction to $N_{C'}-\bdry N_{C'}$ is a diffeomorphism onto $\overline X-C'$, $\psi_{C'}(X\cap N_{C'})=X$ and $\psi_{C'}(C\cap N_{C'})=C$.
Assume that the metric on $X\cap N_{C'}-\bdry N_{C'}$ defined by $d$ and $\psi_{C'}$ is continuously extended to $(X\cap N_{C'})\cup\bdry N_{C'}$.
If we normally modify $\psi:N_C\to\overline X$ then $\bdry N_C$ intersects transversally each face of $N_{C'}$ in $\overline X$, $\psi(N_{C'}\cap N_C)=N_{C'}$ and the metric on $X\cap N_{C'}\cap N_C-\bdry(N_{C'}\cap N_C)$ defined by $d$ and $\psi_{C'}\circ\psi$ is continuously extended to $(X\cap N_{C'}\cap N_C)\cup\bdry(N_{C'}\cap N_C)$.}\vskip1mm
We can say that $d$ becomes continuously extendable to $C$ through $\psi_{C'}\circ\psi$ since $\bdry (N_{C'}\cap N_C)=(\psi_{C'}\circ\psi)^{-1}(C)$.\vskip2mm\noindent
{\it Proof.} We proceed with the proof of the first statement as in the proof of Lemma 1.
We need to define $C'$.
We can replace $(\overline X,C)$ by $(\pi^{-1}(c)\times C,\{c\}\times C)$ for the above $\pi:V\to C$ as the problem depends only on $\pi$ and is local at $C$.
Moreover, we assume $(\overline X,\overline X-X,C)=(B(n-k)\times R^k,\{0\}\times B(l)\times R^k,\{0\}\times R^k),\ 0<k<n-1,\ l<n-k$, for simplicity of notation, where $B(i)$ is the closed unit disk in $R^i$ with center 0.
For each $t\in(0,\,1]$ and each $y\in R^k$, set
$$\delta_0(t,y)=\sup\{d((x,y),(s x,y)):x\in B(n-k)-\{0\}\times B(l),\,s\in(0,\,1],\,|x|=t\}.$$
Note that $\sup_{s\in(0,\,1]}d((x,y),(s x,y))$ is the supremum of the metrics of $(x,y)$ and points in the half-open segment joining $(x,y)$ and $(0,y)$.
Fix $y$.
Then $\delta_0(t,y)$ is a definable continuous function on $(0,\,1]$ and $\lim_{t\to0}\delta_0(t,y)=0$ because after a hole-blow-up of $(B(n-k)-\{0\}\times B(l))\times\{y\}$ with center $(0,y)$, $d|_{(B(n-k)-\{0\}\times B(l))\times\{y\}}$ becomes continuously extendable to the boundary in $B(n-k)\times\{y\}$ of the manifold with boundary of the hole-blow-up as shown near the beginning of Definition 1 of hole-blows-up.
Let $\delta_0$ denote the extension of $\delta_0$ to $[0,\,1]\times R^k$ defined to be 0 on $\{0\}\times R^k$, which is definable but not necessarily continuous at $\{0\}\times R^k$.
Let $C'_0$ be the closure of the subset of $C$ where $\delta_0$ is not continuous.
Then $C'_0$ is definable and of dimension$\,<k$ by the same reason as in the proof of Lemma 1.
As there, we need functions other than $\delta_0$.
For $i=1,\ldots,k$, we define $\delta_i$ by the half-open segments joining $(x,y)$ and $(0,y+|x|e_i)$ likewise $\delta_0$ by the half-open segments joining $(x,y)$ and $(0,y)$, where $e_i$ denotes the $i$-th unit vector in $R^k$.
We also define $C'_i$ by $\delta_i$ in the same way, which is definable and of dimension$\,<k$.
Set $C'=\cup C'_i$.
We show that $C'$ satisfies the condition in Lemma 3.
This is equivalent to the following statement:\vskip1mm\noindent
{\it Set $N=\{(x,y)\in R^{n-k}\times R^k:|x|\ge1\}$, and let $d$ be a definable metric on $\Int N-\{0\}\times R^{l+k}$.
Assume two conditions.
The first is that for each $y\in R^k$, the restriction of $d$ to $(R^{n-k}-\{0\}\times R^l)\times\{y\}\cap\Int N$ is continuously extended to $((N-\{0\}\times R^l)\cup\partial N)\cap R^{n-k}\times\{y\}$.
We define a definable function $\delta_0$ on $[0,\,1]\times R^k$ by the half-open segments joining $(x,y)$ and $(x/|x|,y)$ likewise the above $\delta_0$ by the half-open segments joining $(x,y)$ and $(0,y)$.
Next, for $i=1,\ldots,k$ also, we define a definable function $\delta_i$ on $[0,\,1]\times R^k$ by the half-open segments joining $(x,y)$ and $(x/|x|,y+(|x|-1)e_i)$.
Then the second condition is that $\delta_0,\ldots,\delta_k$ are continuous.
Under these condition, $d$ is continuously extended to $(N-\{0\}\times R^{l+k})\cup\partial N$.}\par
We can prove this statement in the same way as in the proof of Lemma 1.
We omit the details.
The second statement of Lemma 3 is clear by the above proof and the following obvious fact:\par
Let $Y,\,Z_1$ and $Z_2$ be compact $C^1$ manifolds possibly with corners, satisfying the face assumption and such that $Z_1\subset Y,\ Z_2\subset Y,\ \partial Z_1\subset\partial Y$ and $Z_1$ intersects transversally faces of $Y$ and $Z_2$ in $Y$.
Let $\pi$ be a definable $C^1$ submersive retraction onto $Z_1$ of a neighborhood $V$ of $Z_1$ in $Y$.
Then for sufficiently small $\epsilon>0\in R$, the set $\{y\in V:|y-\pi(y)|=\epsilon\}$ is a $C^1$ manifold possibly with corners satisfying the face assumption and it intersects transversally each face of $Y$ and $Z_2$ in $Y$.\qed\vskip2mm
\section{Proof of Theorems}
{\it Proof of Theorem 1.} The uniqueness and the ``only if\," part are obvious.
Hence we only prove that $X$ is compactifiable under the condition in Theorem 1.
Note that $d$ is bounded for the following reason:\par
Assume the contrary.
Then by the curve selection lemma, there are definable continuous curves $\gamma,\gamma':(0,\,1]\to X$ such that $d(\gamma(t),\gamma'(t))\to\infty$ as $t\to0$.
We have
$$d(\gamma(t),\gamma'(t))\le d(\gamma(t),\gamma(1/2))+d(\gamma'(t),\gamma'(1/2))+d(\gamma(1/2),\gamma'(1/2)).$$
By the condition in Theorem 1, $d(\gamma(t),\gamma(1/2))$ and $d(\gamma'(t),\gamma'(1/2)),\ t\in(0,\,1/2]$, are bounded.
Hence $d(\gamma(t),\gamma'(t)),\ t\in(0,\,1/2]$, is bounded, which is a contradiction.\par
We can assume $X=(0,\,1]\times[-1,\,1]^{n-1}$ after the following four steps:
Let $X$ be contained and bounded in $R^n$.
Then we first assume that $X$ is locally closed in $R^n$ because we can replace $X$ by $X-(\overline X-X\overline)$ by uniqueness of the compactification and $X-(\overline X-X\overline)$ is locally closed, where the symbol $(\ \ \overline)$ stands for the closure in $R^n$ of a subset $(\ \ )$ of $R^n$.
Secondly, we suppose that $\overline X-X=\{0\}$ by replacing $X$ by the definable set $\{(\phi(x)x,\phi(x))\in R^{n+1}:x\in X\}$ for a definable continuous function $\phi$ on $\overline X$ with zero set $\overline X-X$.
Thirdly, $X=(0,\,1]\times\{x\in X:|x|=1\}$ because by the o-minimal version of Thom's first isotopy lemma, there is a definable homeomorphism $\pi:(0,\,1]\times\{x\in X:|x|=\epsilon\}\to\{x\in X:|x|\le\epsilon\}$ for small $\epsilon>0$ such that $|\pi(t,x)|=\epsilon t$.
Fourthly, $X=(0,\,1]\times[-1,\,1]^{n-1}$ as follows:
By the triangulation theorem of definable sets, we assume that $\{x\in X:|x|=1\}$ is the underlying polyhedron of a finite simplicial complex $K$.
For each simplex $\sigma$ of $K$, we will choose a definable continuous embedding $\pi_\sigma:(0,\,1]\times\sigma\to R^{n_\sigma}$ of compactification of $(0,\,1]\times\sigma$.
Then we can extend $\pi_\sigma$ to a definable bounded continuous map $\pi_\sigma:X\to R^{n_\sigma}$.
Hence the map $X\ni x\to\Pi_{\sigma\in K}\pi_\sigma(x)\in\Pi_{\sigma\in K}R^{n_\sigma}$ is a pseudo-compactification of $X$, which means that the closure of its image is a pseudo-metric space.
We will show that a pseudo-compactification can be modified to a compactification later.
Thus it is sufficient to prove the proof of Theorem 1 in the case $X=(0,\,1]\times\sigma$ and hence $X=(0,\,1]\times[-1,\,1]^{n-1}$.\par
Let $r$ be a small definable increasing continuous function on $(0,\,1]$ satisfying the condition in Lemma 2, i.e., there is no $r$-bad point.
We will construct a definable continuous map $\psi$ from a compact definable $C^1$ manifold $Y$ with corners to $\overline X$ such that $Y\subset\overline X$, $\psi^{-1}(X)$ is dense in $Y$ and the map $\psi|_{\psi^{-1}(X)}:\psi^{-1}(X)\to X$ is a homeomorphism.
We also define a definable metric $d_\psi$ on $\psi^{-1}(X)$ by $d$ and $\psi$ and a definable continuous function $r_\psi$ by $d_\psi$ likewise $r$ by $d$.
Then it suffices to find $Y$ and $\psi$ so that $d_\psi$ is continuously extended to $Y$ or, equivalently, $r_\psi=1$.
Here we can assume that $r(t)=t$, i.e., there is no $t$-bad point, by replacing the germ of $(0,\,1]\times[-1,\,1]^{n-1}$ at $\{0\}\times[-1,\,1]^{n-1}$ through the embedding $(0,\,1]\times[-1,\,1]^{n-1}\ni(t,x)\to(r(t),x)\in(0,\,\infty)\times[-1,\,1]^{n-1}$.\par
Consider the case $n=2$ where we see the essence of the proof.
Let $\{a_1,\ldots,a_k\}$ be the set of bad points of $(X,d)$, and $\psi:[0,\,\infty)\times[-1,\,1]\cap N\to[0,\,\infty)\times[-1,\,1]$ be a hole-blow-up of $X$ with center $\{a_1,\ldots,a_k\}$.
Then the metric $d_\psi$ on $(0,\,\infty)\times[-1,\,1]\cap\Int N$ is continuously extended to $[0,\,\infty)\times[-1,\,1]\cap N$.
Thus a pseudo-compactification is proved in the case $n=2$.\par
Let $n>2$.
Let $C$ denote the subset of $\{0\}\times[-1,\,1]^{n-1}$ of bad points.
Let $\{C_i:i=1,\ldots,k\}$ be a Whitney definable $C^1$ stratification of $C$ such that $\dim C_i\le\dim C_{i'}$ if $i\le i'$.
Let $\{U_i=(|U_i|,\pi_i,\xi_i)\}$ be a controlled definable tube system for $\{C_i\}$.
Let $\epsilon_1,\ldots,\epsilon_k$ be sufficiently small positive numbers such that $\epsilon_1\gg\cdots\gg\epsilon_k$, and set
$$V_i=\{x\in U_i\cap\overline X:\xi_i(x)\le\epsilon_i\}-\cup_{i'<i}\{x\in U_{i'}:\xi_{i'}(x)<\epsilon_{i'}\}.$$
Then for each $i$, $(X-\cup_{i'<i}V_{i'}\overline)$ and $(C_i-\cup_{i'<i}V_{i'}\overline)$ are compact definable $C^1$ manifolds possibly with corners satisfying the face assumption, and the pair $\big((X-\cup_{i'<i}V_{i'}\overline),(C_i-\cup_{i'<i}V_{i'}\overline)\big)$ admits a hole-blow-up $\psi_i:(X-\cup_{i'\le i}V_{i'}\overline)\to(X-\cup_{i'<i}V_{i'}\overline)$.
(See Figure 1.)
Here $\psi_i$ is defined by the definable $C^1$ submersive retraction $\pi_i|_{V_i}:V_i\to(C_i-\cup_{i'<i}V_{i'}\overline)$.
We will normally modify $\psi_i$.
To be precise, we will shrink not only $V_i$ but also $C_i$ to $\hat V_i$ and $\hat C_i$, respectively, so that $\hat V_i$ is defined by $\hat\epsilon_i$ smaller than $\epsilon_i$, and we replace $\pi_i|_{V_i}$ by $\pi_i|_{\hat V_i\cap\pi^{-1}_i(\hat C_i)}$.
In this case also, if $d_{\psi_i}$ restricted to $(X-V_i)\cap\pi^{-1}_i(x)$ is continuously extended to $((X-V_i)\cup\partial V_i)\cap\pi^{-1}_i(x)$ for any $x\in(C_i-\cup_{i'<i}V_{i'}\overline)$ then this property continues to hold after this normal modification of $\psi_i$.\par
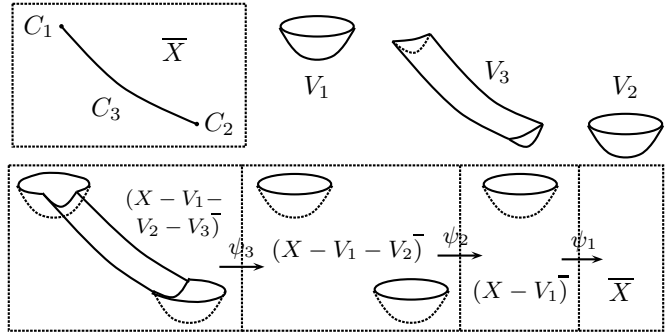
\begin{figure}[!h]
\begin{pspicture}(-2.25,-0.75)(10.5,3.6)
\psline[linestyle=dashed, dash=1pt 0.5pt](0.85,3.6)(0.85,1.7)(3.95,1.7)(3.95,3.6)(0.85,3.6)
\pscurve(1.5,3.3)(2.3,2.5)(3.3,2)\rput(3,3){$\overline X$}\rput(1.2,3.3){$C_1$}\rput(2.1,2.2){$C_3$}\rput(3.6,2){$C_2$}\pscircle*(1.5,3.3){0.03}\pscircle*(3.3,2){0.03}

\psellipse(4.9,3.3)(0.5,0.15)
\pscurve(4.4,3.28)(4.7,2.9)(5.1,2.9)(5.4,3.28)\rput(4.9,2.5){$V_1$}

\psellipse(9,2)(0.5,0.15)
\pscurve(8.5,1.98)(8.8,1.6)(9.2,1.6)(9.5,1.98)\rput(9,2.5){$V_2$}

\pscurve(6.4,3.2)(7.25,2.3)(7.9,2)
\pscurve(5.9,3.1)(6.8,2.15)(7.45,1.8)
\pscurve[linestyle=dashed, dash=1pt 0.5pt](6.4,3.2)(6.2,2.95)(5.9,3.1)
\pscurve(6.4,3.2)(6.2,3.13)(5.9,3.1)
\pscurve(7.9,2)(7.7,1.7)(7.45,1.8)
\pscurve(7.9,2)(7.6,1.83)(7.45,1.8)\rput(7.3,2.7){$V_3$}

\pscurve[linestyle=dashed, dash=1pt 0.5pt](0.91,1.18)(1.2,0.8)(1.6,0.8)(1.88,1.18)
\pscurve(1.25,1.07)(0.93,1.18)(1.18,1.33)(1.4,1.35)(1.62,1.323)(1.87,1.2)(1.7,1.1)\
\pscurve(3.2,-0.07)(3.5,-0.11)(3.68,-0.23)(3.5,-0.34)(3,-0.34)(2.7,-0.25)
\pscurve[linestyle=dashed, dash=1pt 0.5pt](2.72,-0.22)(3,-0.6)(3.4,-0.6)(3.68,-0.22)
\pscurve(1.7,1.1)(2.55,0.3)(3.2,-0.1)
\pscurve(1.25,1.08)(2.15,0.15)(2.75,-0.25)
\pscurve(1.7,1.1)(1.5,0.9)(1.25,1.08)
\pscurve(3.2,-0.05)(3,-0.3)(2.72,-0.22)
\psline[linestyle=dashed, dash=1pt 0.5pt](0.8,1.45)(3.9,1.45)(3.9,-0.75)(0.8,-0.75)(0.8,1.45)
\rput(3,1){\footnotesize$(X-V_1-$}
\rput(3.05,0.65){\footnotesize$V_2-V_3\overline)$}
\rput(3.9,0.35){\small$\psi_3$}\psline{->}(3.6,0.1)(4.2,0.1)

\psellipse(4.6,1.2)(0.5,0.15)
\pscurve[linestyle=dashed, dash=1pt 0.5pt](4.1,1.18)(4.4,0.8)(4.8,0.8)(5.1,1.18)
\psellipse(6.15,-0.2)(0.5,0.15)
\pscurve[linestyle=dashed, dash=1pt 0.5pt](5.65,-0.22)(5.95,-0.6)(6.35,-0.6)(6.65,-0.22)
\psline[linestyle=dashed, dash=1pt 0.5pt](3.9,1.45)(6.8,1.45)(6.8,-0.75)(3.9,-0.75)
\rput(5.3,0.35){\small$(X-V_1-V_2\overline)$}
\rput(6.75,0.5){\small$\psi_2$}\psline{->}(6.5,0.25)(7.1,0.25)

\psellipse(7.6,1.2)(0.5,0.15)
\pscurve[linestyle=dashed, dash=1pt 0.5pt](7.1,1.18)(7.4,0.8)(7.8,0.8)(8.1,1.18)
\psline[linestyle=dashed, dash=1pt 0.5pt](6.8,1.45)(8.38,1.45)(8.38,-0.75)(6.8,-0.75)
\rput(7.6,-0.2){\small$(X-V_1\overline)$}
\rput(8.45,0.45){\small$\psi_1$}\psline{->}(8.15,0.2)(8.7,0.2)

\rput(8.94,-0.2){$\overline X$}
\psline[linestyle=dashed, dash=1pt 0.5pt](8.38,1.45)(9.5,1.45)(9.5,-0.75)(8.38,-0.75)
\end{pspicture}
\caption{Hole-blow-up with general center}
\end{figure}\par
We modify $\{C_i\}$ and $\{U_i\}$ so that $d_{\psi_i}$ on $X-\cup_{i'\le i}V_{i'}$ defined by $d$ and $\psi_i$ is continuously extended to $(X-\cup_{i'\le i}V_{i'})\cup(\partial V_i-\cup_{i'<i}V_{i'})$ for any $\{\epsilon_i\}$.
We proceed by downward induction.
Assume that $d_{\psi_{i+1}},\ldots,d_{\psi_k}$ satisfy this condition for some $i$.
We need to shrink $C_i$ and enlarge $\cup_{i'=1}^{i-1}C_{i'}$ since $d_{\psi_i}$ does not necessarily satisfy the condition.
By Lemma 3 there exists a definable closed subset $C'_i$ of $C_i$ of dimension $<\dim C_i$ such that $d_{\psi_i}$ on 
$X-\cup_{i'\le i}V_{i'}$ is continuously extended to $(X-\cup_{i'\le i}V_{i'})\cup(\partial V_i-\cup_{i'<i}V_{i'})-\psi^{-1}_i(C')$.
Replace $C_i$ by $C_i-C'_i$ and substratify $\{C_i-C'_i,C'_i,\cup_{i'<i}C_{i'}\}$ to a Whitney definable $C^1$ stratification $\{C''_{i'}:i'\le i\}$, although the number of $C_1,\ldots,C_i$ may change but we ignore it.
Then $d_{\psi_i},\ldots,d_{\psi_k}$ satisfy the condition even if $\psi_i,\ldots,\psi_k$ are normally modified and hence when we modify $U_i-\pi^{-1}_i(C'_i),U_{i+1},\ldots,U_k$ for any $\epsilon_i,\ldots,\epsilon_k$.
This property is important.
Indeed, there does not necessarily exist a definable tube system for $\{C''_{i'}:i'<i\}$ together with which $\{U_i-\pi^{-1}_i(C'_i),U_{i'}:i'>i\}$ is controlled.
(See the proof of Lemma I.1.3 in \cite{S1}.)
Hence we need to modify $\{U_i-\pi^{-1}_i(C'_i),U_{i'}:i'>i\}$ so that this is the case.
Choose a controlled definable tube system for $\{C''_{i'}:i'\le i\}\cup\{C_{i'}:i'>i\}$.
Then the induction step is proved.
Thus all of $d_{\psi_i}$ satisfy the condition.\par
Under the last condition on $d_{\psi_i}$, we show that $d_{\psi_1\circ\cdots\circ\psi_k}$ is continuously extended to $(X-\cup_{i=1}^kV_i\overline)$ by induction.
Note $(X-\cup_{i=1}^kV_i\overline)=(X-\cup_{i=1}^kV_i)\cup\bdry\cup_{i=1}^kV_i$.
Assume that $d_{\psi_1\circ\cdots\circ\psi_i}$ is continuously extended to $(X-\cup_{l=1}^iV_i)\cup\bdry\cup_{l=1}^iV_l$ for some $i$.
Then $d_{\psi_1\circ\cdots\circ\psi_{i+1}}$ is continuously extended to $(X-\cup_{l=1}^{i+1}V_l)\cup\bdry\cup_{l=1}^{i+1}V_l$ after some normal modification of $\psi_{i+1}$ by the second statement of Lemma 3.
Thus the induction step and the pseudo-compactification are proved.
Hence Theorem 1 follows from the next lemma, which says that a pseudo-compactifiable definable metric space is compactifiable.\qed\vskip2mm\noindent
{\bf Lemma 4.} {\it Let $X$ be a compact definable set with an equivalence relation $\sim$ such that the set $Y=\{(x_1,x_2)\in X^2:x_1\sim x_2\}$ is definable and closed in $X^2$.
Then the quotient space $\tilde X$ admits a definable set structure such that the projection $q:X\to\tilde X$ is definable.}\vskip2mm\noindent
{\it Proof.} Let an equivalence relation on a subset $X_0$ of $X$ be always the restriction of the equivalence relation to $X_0$, let $\widetilde{X_0}$ denote the quotient space of $X_0$, and let $q_{X_0}:X_0\to\widetilde{X_0}$ denote the projection.
We call $X_0$ closed under the equivalence relation if any point of $X_0$ is not equivalent to any point of $X-X_0$.
Note that given a definable closed subset $C$ of $X$ and the $q_C:C\to\tilde C$ as required in Lemma 4 (i.e., the condition in Lemma 4 holds for $\tilde C$ and $q_C$), $q_C$ is extended to a definable continuous map $g_C:X\to R^n$ for some $n$ so that $g_C|_{X-C}$ is a homeomorphism onto $g(X)-g(C)$.
(Here we enlarge the ambient Euclidean space of $\tilde C$, and we do not mention the enlarging in the argument below.)
The reason is the following:\par
If we ignore the condition that $g_C|_{X-C}$ is a homeomorphism, we can extend $q_C$ to a definable continuous map $g_1:X\to R^n$.
We need to modify $g_1$ so that the condition holds.
Let $\phi$ be a definable continuous function on $X$ with zero set $C$.
Let $\id$ denote the identity map of $X$, and set $g_2(x)=(\phi(x)\id(x),\phi(x))$ for $x\in X$.
Then $g_2$ is a definable continuous map from $X$ to some $R^{n'}$, $g_2(C)=\{0\}$, $g_2(C)\cap g_2(X-C)=\emptyset$ and $g_2|_{X-C}$ is a homeomorphism onto $g(X)-g(C)$.
Hence $g_C$, defined to be $(g_1,g_2)$, satisfies the requirements.\par
Assume that Lemma 4 is already proved for $X$ of dimension $<k$.
Then we need to prove Lemma 4 for $X$ of dimension $k$.
We prove the following more general statement:\vskip1mm\noindent
{\it Given definable closed subsets $A_1\subset A_2$ of $X$ and the map $q_{A_1}:A_1\to\widetilde{A_1}$ as required, $q_{A_1}$ is extended to $q_{A_2}:A_2\to\widetilde{A_2}$ as required.}\par
We prove this statement by induction on $\dim A_2$.
For this, it suffices to prove the statement in the case $A_2=X$ by the above note.
Thus we assume that $A_2=X$, and the induction hypothesis is that the statement holds for $X$ of dimension $<k$.
The proof in the case $\dim A_1=k$ is essentially the same as in the case $\dim A_1<k$.
Hence we prove the statement in the case $\dim A_1<k$.
Let $p:Y\to X$ denote the restriction to $Y$ of the projection of $X^2$ onto the first factor.
Let $\{Y_j\}$ and $\{X_i\}$ be respective Whitney definable $C^1$ stratifications of $Y$ and $X$ such that $A_1$ is the union of some $X_i$'s, the set of definably $C^1$ smooth points of $p^{-1}(X_i)$ of dimension $=\dim p^{-1}(X_i)$ is the union of some $Y_j$'s for each $i$, each connected component of $X_i$ is contractible and of dimension $i$ and the restriction of $p$ to each connected component of $Y_j$ is a surjective $C^1$ submersion onto some connected component of $X_i$ (see Lemma I.2.6 in \cite{S1} for the proof of existence of such stratifications).
We write $p:\{Y_j\}\to\{X_i\}$.
Then for each connected component of $X_{i0}$ of $X_i$, we have a definable homeomorphism $h:X_{i0}\times p^{-1}(x_{i0})\to p^{-1}(X_{i0}),\ x_{i0}\in X_{i0}$, such that $p\circ h:X_{i0}\times p^{-1}(x_{i0})\to X_{i0}$ is the projection and for each $Y_j$, $h|_{X_{i0}\times(p^{-1}(x_{i0})\cap Y_j)}$ is a $C^1$ diffeomorphism onto $p^{-1}(X_{i0})\cap Y_j$ by the o-minimal version of Thom's first isotopy lemma, since $X_{i0}$ is contractible.
Namely $p$ is trivial over $X_{i0}$.
We will shrink $X_{i0}$.
Then we loose the property that it is contractible but keep the triviality property.
Note that if $X_{i0}$ is closed under the equivalence relation and if $\dim X_{i0}<\dim p^{-1}(X_{i0})$ then $p^{-1}(X_{i0})\subset X_{i0}\times X_{i0}$, $\widetilde{X_{i0}}$ admits a definable set structure such that $q_{X_{i0}}:X_{i0}\to\widetilde{X_{i0}}$ is definable, $\widetilde{X_{i0}}$ is of dimension $<i$ and there is a homeomorphism from $X_{i0}$ to $\widetilde{X_{i0}}\times p^{-1}(x),\ x\in X_{i0}$, which preserves the equivalence relations, where the equivalence relation on $\widetilde{X_{i0}}\times p^{-1}(x)$ is defined so that elements $(x_1,y_1)$ and $(x_2,y_2)$ of $\widetilde{X_{i0}}\times p^{-1}(x)$ are equivalent if and only if $x_1=x_2$.\par
By the induction hypothesis we have $q_{X-X_k}:X-X_k\to(X-X_k\tilde)$ which is an extension of $q_{A_1}$, where $(X-X_k\tilde)$ means the quotient space of $X-X_k$, and we need to extend $q_{X-X_k}:X-X_k\to(X-X_k\tilde)$ to $q:X\to\tilde X$.
After each step of the construction below of $q_A$, $A\subset X$, we need to define its extension $g_A:X\to R^n$.
However, we do not mention it.
First we reduce the problem to the case $\overline{X_k}=X$ and $X_k$ is connected and closed under the equivalence relation.
This is possible for the following reason:\par
Let $X'_k$ denote the subset of $X_k$ consisting of $x$ which is equivalent to some $x'$ of $X-X_k$.
Then $X'_k$ is definable and the map $X'_k\ni x\to q_{X-X_k}(x')\in(X-X_k\tilde)$ is well defined and definable since $Y$ is definable.
Moreover, its extension to $X'_k\cup(X-X_k)$, defined to be $q_{X-X_k}$ on $X-X_k$, is continuous.
Indeed, let $C$ be a definable closed subset of $(X-X_k\tilde)$.
Then $q^{-1}_{X-X_k}(C)$ is definable and closed in $X-X_k$ since $q_{X-X_k}:X-X_k\to(X-X_k\tilde)$ is definable and continuous, hence $Y\cap(X\times q^{-1}_{X-X_k}(C))$ is definable and closed in $Y$, $p(Y\cap(X\times q^{-1}_{X-X_k}(C))$ is closed in $X$, and this set equals the subset of $X'_k\cup(X-X_k)$ consisting of points equivalent to some points in $q^{-1}_{X-X_k}(C)$.
(Here we use the fact that the image of a definable continuous map from a definable compact set into a Euclidean space is compact, although a definable compact set is not necessarily compact in the topological sense.)
By the same reason, $X'_k$ is closed in $X_k$.
Hence $X_k-X'_k$ is a definable $C^1$ manifold of dimension $k$ and closed under the equivalence relation, $q_{X-X_k}$ is extended to $q_{X'_k\cup(X-X_k)}:X'_k\cup(X-X_k)\to(X'_k\cup(X-X_k)\tilde)$ as required, and we have $q_{\overline{X_{k0}}-X_{k0}}:\overline{X_{k0}}-X_{k0}\to(\overline{X_{k0}}-X_{k0}\tilde)$, where $X_{k0}$ is a connected component of $X_k-X'_k$.
Therefore, we can replace $(X,X_k)$ by $(\overline{X_{k0}},X_{k0})$.
(Here there is the case where some point of $X_{k0}$ is equivalent to some point of $X_k-X_{k0}$.
Then we construct an extension $q_{X_{k0}\cup(X-X_k)}$ first, next for another connected component $X_{k1}$ of $X_k$, we extend it further to $q_{X_{k0}\cup(X-X_k)\cup X'_{k1}}$, then to $q_{X_{k0}\cup(X-X_k)\cup X_{k1}}$, and so on, where $X'_{k1}$ is the subset of $X_{k1}$ of points equivalent to some points in $X_{k0}$.)
Consequently, the above reduction is possible.\par
Set $\{Y_{k j}:j\}=\{Y_j:p(Y_j)=X_k\}$.
We suppose $\dim Y_{k j}=j$ by gathering $Y_{k j}$ of the same dimension.
Next we can assume that $\{Y_{k j}:j\}$ consists of one stratum, say, $Y_{k k_1}$ for the following reason:\par
Assume that the set is not a singleton.
We come back to the definition of $p:\{Y_j\}\to\{X_i\}$.
There we have chosen $\{Y_{k j}\}$ so that for each $x\in X_k$, the union of $Y_{k j}\cap p^{-1}(x)$ of the maximal dimension is the set of definably $C^1$ smooth points of $Y\cap p^{-1}(x)$ of dimension $=\dim Y\cap p^{-1}(x)$.
Let $Z_k$ denote the union.
Then $X_k-p_2(Z_k)$ is nonempty, closed in $X_k$ and of dimension smaller than $k$, where $p_2$ is the restriction to $Y$ of the projection $X\times X\to X$ onto the second factor because the $p_2(p^{-1}(x))$ is the subset of $X_k$ of points equivalent to $x$, the family $\{p_2(p^{-1}(x)):x\in X_k\}$ is a definable $C^0$ foliation of $X_k$ and $p_2(Z_k)$ is the union of the sets of definably $C^1$ smooth points of $p_2(p^{-1}(x))$ of dimension $=\dim p_2(p^{-1}(x))$ for $x\in X_k$.
Hence by the induction hypothesis, we can extend $q_{X-X_k}$ to $q_{X-p_2(Z_k)}:X-p_2(Z_k)\to(X-p_2(Z_k)\tilde)$.
Then we can extend it further to $q:X\to\tilde X$ since any point of $p_2(Z_k)$ is equivalent to some point of $X_k-p_2(Z_k)$.
Thus Lemma 4 is proved in the case where $\{Y_{k j}:j\}$ is not a singleton.\par
It only remains to prove the case $\{Y_{k j}:j\}=\{Y_{k k_1}\}$.
By the triviality property, $Y_{k k_1}$ is definably $C^1$ diffeomorphic to $X_k\times p^{-1}(x_k),\ x_k\in X_k$, and the family $\{p_2(p^{-1}(x)):x\in X_k\}$ is a definable product $C^1$ foliation of $X_k$.
Let $h:p_2(p^{-1}(x_k))\times F\to X_k$ be a definable $C^1$ diffeomorphism for fixed $x_k\in X_k$ and some definable $C^1$ manifold $F$ such that $\{h\big(p_2(p^{-1}(x_k))\times\{y\}\big):y\in F\}=\{p_2(p^{-1}(x)):x\in X_k\}$.
There are two case: $\dim p^{-1}(x_k)\,(=\dim p_2(p^{-1}(x_k)))>0$, i.e., the foliation is of positive dimension or $\dim p^{-1}(x_k)=0$.
Consider the first case.
The set $h(\{x_k\}\times F)$ is definable, closed in $X_k$ and of dimension $<k$, the equivalence relation restricted to $h(\{x_k\}\times F)$ is the equality equivalence relation, and hence $(h(\{x_k\}\times F)\tilde)=h(\{x_k\}\times F)$.
Therefore, we have an extension $q_{(X-X_k)\cup h(\{x_k\}\times F)}:(X-X_k)\cup h(\{x_k\}\times F)\to\big((X-X_k)\cup h(\{x_k\}\times F)\tilde{\big)}\,(=(X-X_k\tilde)\cup h(\{x_k\}\times F))$, which is equal to $q_{X-X_k}$ on $X-X_k$ and to the identity map on $h(\{x_k\}\times F)$.
We need to extend it further to $q:X\to\tilde X$.
Each leaf $p_2(p^{-1}(x))$ contains one and only one point $x'$ in $h(p_2(y_k)\times F)$ and the correspondence from $x$ to $x'$ is a definable continuous map from $X_k$ to $p_2\circ h(\{x_k\}\times F)$.
Hence we can extend $q_{(X-X_k)\cup h(\{x_k\}\times F)}$ to $q:X\to((X-X_k)\cup h(\{x_k\}\times F)\tilde{\big)}$ by setting $q(x)=q_{(X-X_k)\cup h(\{x_k\}\times F)}(x')$.\par
If $\dim p^{-1}(x_k)=0$ then $\widetilde{X_k}=X_k$.
Hence the map $g_{X-X_k}:X\to g_{X-X_k}(X)$ satisfies the conditions on $q$.
In any case we have $q$, which is an extension of $q_{A_1}$ since $A_1$ is the union of some $X_i$'s.
Thus the proof is complete.\qed\vskip2mm\noindent
{\it Proof of Theorem 2.}
For simplicity of presentation we proceed with the o-minimal structure $\R_{\rm an}$, and later show the semialgebraic case.
Set $B=\{x\in\R^2:|x|\le 1\}$, and regard $\R^2$ as $\C$.
We will give a subanalytic metric $d$ on $B$.
Before this, we set $X=(0,\,1]\times\R$, and we will define a semialgebraic metric $d_X$ on $X$ uniform on the factor $\R$, i.e., $d_X(x,x')=d_X(x+a,x'+a)$ for $x,x'\in X$ and $a\in\{0\}\times\R$ and so that the following condition holds:\\
$(*)$ {\it For $b,b'\in B-\{0\}$, $(x_1,x_2)\in f^{-1}(b)$ and $(x'_1,x'_2)\in f^{-1}(b')$,
\begin{align*}d_X(f^{-1}(b),f^{-1}(b'))=d_X((x_1,x_2),(x'_1,x'_2))\ \text{ if }\,|x_2-x'_2|\le1,\\
d_X(f^{-1}(b),f^{-1}(b'))\le d_X((x_1,x_2),(x'_1,x'_2))\ \text{ otherwise},\ \ \qquad\end{align*}}\noindent
where $f:X\to B$ is an analytic map defined by $f(x_1,x_2)=x_1\exp(\pi i x_2)$ for $(x_1,x_2)\in X$.
(Here $\pi$ and $i$ in $\exp(\pi i x)$ are the ratio of the circumference of a circle and the imaginary unit respectively.)
Let $d|_{B-\{0\}}$ be induced by $d_X$ through $f$, to be precise, $d(b,b')=d_X(f^{-1}(b),f^{-1}(b'))$ for $b,b'\in B-\{0\}$.
Then $d|_{B-\{0\}}$ is a well-defined semialgebraic metric because $d_X$ satisfies Condition $(*)$ and is uniform on $\R$ and $f((0,\,1]\times[0,\,2))=B-\{0\}$.
We extend $d|_{B-\{0\}}$ to $d$ on $B$ by setting $d(0,b)=|b|/2$ for $b\in B$.
Then we require $d$ to be a subanalytic metric.
For this, the condition that $d_X$ satisfies is
$$\lim_{b\to0}d_X(f^{-1}(b),f^{-1}(b'))=|b'|/2\ \ \text{for }b,b'\in B-\{0\}.\leqno{(**)}$$
For the time being, we assume that there exists $d_X$ satisfying Conditions $(*)$ and $(**)$.
This $(B,d)$ is one of the subanalytic metric space.\par
We want to define another subanalytic metric $\tilde d$ on the same $B$ using a new semialgebraic metric $\tilde d_X$ on $X$ by a similar method.
Let $\phi$ be a semialgebraic function on $(0,\,1]$ such that $\phi$ is strictly increasing, $\lim_{t\to0}\phi(t)=\infty$ and $\phi(1)=0$.
For example, $\phi$ is defined by $\phi(t)=-1/t+1$.
We define a semialgebraic metric $\tilde d_X$ on $X$ by $\tilde d_X((x_1,x_2),(x'_1,x'_2))=d_X\big((x_1,x_2+\phi(x_1)),(x'_1,x'_2+\phi(x'_1))\big)$ for $(x_1,x_2),(x'_1,x'_2)\in X$ and then a function $\tilde d$ on $B^2$ by $\tilde d(0,b)=|b|/2$ for $b\in B$ and $\tilde d(b,b')=\tilde d_X(f^{-1}(b),f^{-1}(b'))$ for $b,b'\in B-\{0\}$.
Then $\tilde d$ is a subanalytic metric by the same reason as for $d$ because $\tilde d_X$ also satisfies the following conditions $\widetilde{(*)}$ and $\widetilde{(**)}$ and is uniform on $\R$ and $f\{(x_1,x_2+\phi(x_1)):(x_1,x_2)\in(0,\,1]\times[0,\,2)\}=B-\{0\}$.\vskip1mm\noindent
$\widetilde{(*)}$ {\it For the same $b,b',(x_1,x_2)$ and $(x'_1,x'_2)$ as in Condition $(*)$,
\begin{align*}\tilde d_X(f^{-1}(b),f^{-1}(b'))=\tilde d_X((x_1,x_2),(x'_1,x'_2))\ \text{ if }\,|x_2-x'_2+\phi(x_1)-\phi(x'_1)|\le1,\\
\tilde d_X(f^{-1}(b),f^{-1}(b'))\le\tilde d_X((x_1,x_2),(x'_1,x'_2))\ \text{ otherwise}.\ \ \qquad\end{align*}}\vskip-6mm
$$\lim_{b\to0}\tilde d_X(f^{-1}(b),f^{-1}(b'))=|b'|/2\ \ \text{for }b,b'\in B-\{0\}.\leqno{\widetilde{(**)}}$$
Here Condition $\widetilde{(*)}$ holds obviously, and Condition $\widetilde{(**)}$ does because its left side equals the first term of the following equalities:
\begin{align*}\lim_{b\to0}\inf_{(x_1,x_2)\in f^{-1}(b)\atop{(x'_1,x'_2)\in f^{-1}(b')}}\tilde d_X((x_1,x_2),(x'_1,x'_2))=\hskip20mm\\
\lim_{b\to0}\inf_{(x_1,x_2)\in f^{-1}(b)\atop{(x'_1,x'_2)\in f^{-1}(b')}}d_X((x_1,x_2+\phi(x_1)),(x'_1,x'_2+\phi(x'_1)))=|b'|/2.\end{align*}
Let $\tilde B$ denote the topological space $B$ with the metric $\tilde d$, and $B$ mean $B$ with the metric $d$.
In the same way, $X$ and $\tilde X$ stand for the semialgebraic metric spaces $(X,d_X)$ and $(X,\tilde d_X)$.\par
Next we define an orientation-preserving isometry $\Phi_{b_0}:B\to\tilde B$ such that $\Phi_{b_0}(1,0)=b_0$ for each $b_0\in\partial B$.
Let $x_0\in\R$ be such that $f(1,x_0)=b_0$.
Set $\Phi_X(x_1,x_2)=(x_1,x_0+\phi(x_1)+x_2)$ for $(x_1,x_2)\in X$.
Then $\Phi_X$ is a semialgebraic orientation-preserving isometry from $\tilde X$ to $X$ because
\begin{align*}d_X(\Phi_X(x_1,x_2),\Phi_X(x'_1,x'_2))=d_X((x_1,x_0+\phi(x_1)+x_2),(x'_1,x_0+\phi(x'_1)+x'_2))=\qquad\\
d_X((x_1,\phi(x_1)+x_2),(x'_1,+\phi(x'_1)+x'_2))=\tilde d_X((x_1,x_2),(x'_1,x'_2))\text{ for }(x_1,x_2),(x'_1,x'_2)\in X.\end{align*}
Note that $\Phi_X$ carries the level sets of the function $x\to\tilde d_X(x,x_0)$ in Figure 4 to those of the $x\to d_X(x,x_0)$ in Figure 2.
Define a map $\Phi_{b_0}:\tilde B\to B$ by $\Phi_{b_0}(b)=f\circ\Phi_X(f^{-1}(b))$.
Then $\Phi_{b_0}$ is well defined and an isometry from $\tilde B$ to $B$ because
\begin{align*}f\circ\Phi_X(x_1,x_2+2)=f(x_1,x_0+\phi(x_1)+x_2+2)=x_1\exp(\pi i(x_0+\phi(x_1)+x_2+2))=\ \\
x_1\exp(\pi i(x_0+\phi(x_1)+x_2))=f(x_1,x_0+\phi(x_1)+x_2)=f\circ\Phi_X(x_1,x_2)\text{ for }(x_1,x_2)\in X,\end{align*}\vskip-6mm
\begin{align*}d(\Phi_{b_0}(b),\Phi_{b_0}(b'))=d_X\big(f^{-1}(\Phi_{b_0}(b)),f^{-1}(\Phi_{b_0}(b'))\big)=d_X\big(\Phi_X(f^{-1}(b)),\Phi_X(f^{-1}(b'))\big)\\
=\tilde d_X(f^{-1}(b),f^{-1}(b'))=\tilde d(b,b')\ \ \text{ for }b,b'\in B-\{0\}\hskip16mm\text{and}
\end{align*}\vskip-4mm
$$\tilde d(0,b')=\lim_{b\to0}\tilde d(b,b')=\lim_{b\to0}d(\Phi_{b_0}(b),\Phi_{b_0}(b'))=d(0,\Phi_{b_0}(b'))\ \text{ for the same }b,b'.$$
Moreover, $\Phi_{b_0}(1,0)=b_0$.
Indeed,
$$\Phi_{b_0}(1,0)=f\circ\Phi_X(f^{-1}(1,0))=f\circ\Phi_X(1,0)=f(1,x_0)=b_0.$$
(The $\Phi_{b_0}$ carries the level sets of the function $b\to\tilde d(b,b_0)$ in Figure 4 to those of the $b\to d(b,b_0)$ in Figure 3.)\par
Now we choose the semialgebraic metric $d_X$ so that it is uniform on $\R$ and satisfies Conditions $(*)$ and $(**)$.
For $x=(x_1,x_2),x'=(x'_1,x'_2)\in X$ with $x'_1\le x_1$, set
$$d_X(x,x')=\left\{\begin{array}{l}
\max\{x_1-x'_1,|x_2-x'_2|\}\ \text{ if }x_1-x'_1\le x_1/2\text{ and }|x_2-x'_2|\le x_1/2\\
x_1/2\hskip32mm\text{otherwise}.\end{array}\right.\nonumber$$
If $x'_1>x_1$, we define $d_X(x,x')$ to be $d_X(x',x)$.
(See Figure 2.)
Then we easily see that $d_X$ is a metric, it is clearly semialgebraic and uniform on $\R$, Condition $(**)$ is obvious, and Condition $(*)$ holds for the following reason:
By Figure 2 we see
\begin{align*}d_X(f^{-1}(f(x)),f^{-1}(f(x')))=d_X((x_1,x_2),f^{-1}(f(x')))=\quad\ \ \,\\
\min_{n\in\Z}d_X((x_1,x_2),(x'_1,x'_2+2n))=d_X((x_1,x_2),(x'_1,x'_2+2n_0))\end{align*}
for $n_0\in\Z$ such that $(x_1,x'_2+2n_0)$ is nearest to the segment joining $(0,x_2)$ and $(1,x_2)$.
There are two cases: two points $(x_1,x'_2+2n_0)$ and $(x_1,x'_2+2n_0+2)$ are nearest, or only one $(x_1,x'_2+2n_0)$ is so.
In the first case, $x'_2+2n_0+1=x_2$ hence $x_2-x'_2-2n_0=1$ and Condition $(*)$ holds.
In the second case, $x'_2+2n_0-1<x_2<x'_2+2n_0+1$, $|x_2-x'_2-2n_0|<1$ and Condition $(*)$ holds.\par
\begin{figure}[!h]
\begin{pspicture}(-0.85,1.3)(10.5,5)
\rput(3,3.62){\tiny$x_0$}
\pscircle*(3,3.5){0.03}
\rput(4.5,3.3){\tiny$2x_0$}
\rput(4.55,2.9){\tiny slope}
\rput(4.55,2.7){\tiny $=\!-1$}
\pscircle*(4,3.5){0.03}
\rput(4.6,3.9){\tiny$3x_0/2$}
\pscircle*(3.48,3.5){0.03}
\psline[linestyle=dashed, dash=2pt 1pt](2,2)(2,5)
\pscustom[linewidth=0.1pt]{\psline[linestyle=dashed, dash=0.1pt 0.1pt](2,2)(2,5)\gsave
\psline[linewidth=1pt](3,5)(3,4.5)(2.5,4)(2.5,3)(3,2.5)(3,2)
\fill[linewidth=0.1pt,fillstyle=hlines]\grestore}
\psline(2,5)(3,5)
\psline(2,2)(3,2)
\psline(3,5)(3,4.5)(2.5,4)(2.5,3)(3,2.5)(3,2)
\psline(3,4.5)(3.48,4.02)(3.48,2.98)(3,2.5)
\psline(2.7,3.8)(2.7,3.2)(3,2.9)(3.3,3.2)(3.3,3.8)(3,4.1)(2.7,3.8)
\psline[linestyle=dotted, dash=1pt 0.1pt](2.5,3)(4.15,4.65)
\psline[linestyle=dotted, dash=1pt 0.1pt](2.5,4)(4.15,2.35)
\psline[linestyle=dotted, dash=1pt 0.1pt](3,2.5)(4.15,2.5)
\psline[linestyle=dotted, dash=1pt 0.1pt](3,4.5)(4.15,4.5)
\psline(3.2,5)(3.2,4.5)(3.62,4.12)(3.62,2.88)(3.2,2.5)(3.2,2)
\psline(3.4,5)(3.4,4.5)(3.74,4.2)(3.74,2.78)(3.4,2.5)(3.4,2)
\rput(3.1,1.3){\footnotesize$x_0\in(0,\,2/3)\!\times\!\{0\}$}
\psline(4,5)(4,2)
\psline(4.15,5)(4.15,2)
\rput(1.9,1.65){\scriptsize$\{0\}\!\times\![-1,1]$}
\rput(4.25,1.65){\scriptsize$\{1\}\!\times\![-1,1]$}
\rput(1.6,3.93){\scriptsize{value}}
\rput(1.55,3.7){\scriptsize=}
\rput(1.58,3.5){\scriptsize$|x_0|/2$}
\rput(1.6,2.6){\tiny{slope}}
\rput(1.6,2.4){\tiny$=1$}
\pscurve[linestyle=dashed, dash=1pt 1pt]{->}(1.7,3.3)(1.85,3.2)(2.25,3.2)
\pscurve[linestyle=dashed, dash=1pt 1pt]{->}(4.45,3.43)(4.25,3.5)(4.05,3.5)
\psline[linestyle=dashed, dash=1pt 1pt]{->}(4.5,3.8)(3.55,3.5)
\pscurve[linestyle=dashed, dash=1pt 1pt]{->}(1.9,2.5)(2.5,2.65)(2.65,3.1)

\pscurve[linestyle=dashed, dash=1pt 1pt]{->}(4.5,2.8)(4.25,2.85)(4.05,2.7)
\psline[linestyle=dashed, dash=2pt 1pt](6.2,5)(6.2,2)
\pscustom[linewidth=0.1pt]{\psline[linestyle=dashed, dash=0.1pt 0.1pt](6.2,2)(6.2,5)\gsave
\psline[linewidth=1pt](7.7,5)(7.7,4.6)(7.05,3.95)(7.05,3.04)(7.7,2.4)(7.7,2)
\fill[linewidth=0.1pt,fillstyle=hlines]\grestore}
\psline(6.2,5)(7.7,5)
\psline(6.2,2)(7.7,2)
\psline(7.7,4.6)(8.17,4.13)
\psline(8.17,2.87)(7.7,2.4)
\rput(7.7,3.62){\footnotesize$x_0$}
\psline(7.7,5)(7.7,4.6)(7.05,3.95)(7.05,3.04)(7.7,2.4)(7.7,2)
\rput(7.25,1.3){\footnotesize$x_0\in[2/3,\,1)\!\times\!\{0\}$}
\pscircle*(7.7,3.5){0.04}
\psline(7.3,3.9)(7.3,3.1)(7.7,2.7)(8.1,3.1)(8.1,3.9)(7.7,4.3)(7.3,3.9)
\psline[linestyle=dotted, dash=1pt 0.1pt](8.17,3.97)(7.7,3.5)
\psline[linestyle=dotted, dash=1pt 0.1pt](7.7,3.5)(8.17,3.03)
\psline(7.9,5)(7.9,4.6)(8.17,4.33)
\psline(7.9,2)(7.9,2.4)(8.17,2.67)
\psline(8.1,5)(8.1,4.6)(8.17,4.53)
\psline(8.1,2)(8.1,2.4)(8.17,2.47)
\rput(6.2,1.65){\scriptsize$\{0\}\!\times\![-1,1]$}
\rput(8.17,1.65){\scriptsize$\{1\}\!\times\![-1,1]$}
\rput(5.8,3.93){\scriptsize{value}}
\rput(5.75,3.7){\scriptsize=}
\rput(5.77,3.5){\scriptsize$|x_0|/2$}
\pscurve[linestyle=dashed, dash=1pt 1pt]{->}(5.9,3.3)(6.1,3.2)(6.5,3.2)

\psline[linestyle=dashed, dash=2pt 1pt](10,5)(10,2)
\pscustom[linewidth=0.1pt]{\psline[linestyle=dashed, dash=0.1pt 0.1pt](10,2)(10,5)\gsave
\psline[linewidth=1pt](12,5)(12,4.5)(11,4)(11,3)(12,2.5)(12,2)
\fill[linewidth=0.1pt,fillstyle=hlines]\grestore}
\psline(10,5)(12,5)
\psline(10,2)(12,2)
\psline(12,5)(12,4.5)(11,4)(11,3)(12,2.5)(12,2)
\psline(12,4.15)(11.35,3.8)(11.35,3.15)(12,2.85)
\psline(12,3.85)(11.7,3.7)(11.7,3.3)(12,3.15)
\rput(9.9,1.65){\scriptsize$\{0\}\!\times\![-1,1]$}
\rput(12.1,1.65){\scriptsize$\{1\}\!\times\![-1,1]$}
\rput(12,3.62){\footnotesize$x_0$}
\rput(11.1,1.3){\footnotesize$x_0=(1,0)$}
\pscircle*(12,3.5){0.04}
\rput(9.6,3.9){\scriptsize{value}}
\rput(9.6,3.6){\scriptsize{=$1/2$}}
\pscurve[linestyle=dashed, dash=1pt 1pt]{->}(9.7,3.4)(9.85,3.3)(10.4,3.3)

\end{pspicture}
\caption{Level sets of the function $x\to d_X(x,x_0)$}
\begin{quote}{\small Each level curve of $x\to d_X(x,x_0)$ is a union of vertical segments and segments with slope $\pm1$ if the value is not $|x_0|/2$.
If $d_X(x,x_0)=x_0/2$, the level set is a union of such three segments and a two-dimensional domain.}\end{quote}
\end{figure}
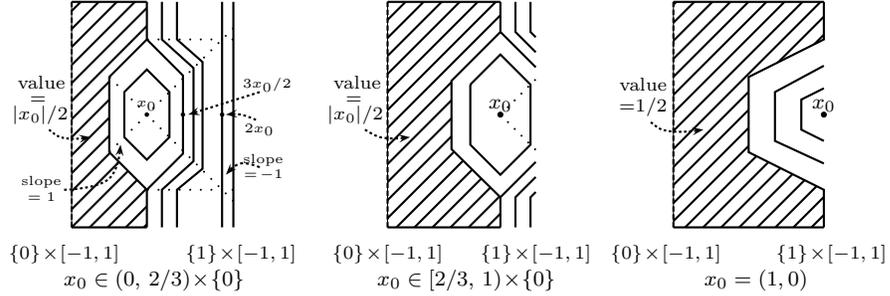

\begin{figure}[!h]
\begin{pspicture}(-0.55,1.6)(10.5,5)
\rput(3.4,3.5){\small0}
\pscircle(3.4,3.5){1.5}
\pscircle(3.4,3.5){1}
\pscircle(3.4,3.5){0.5}
\rput(3.4,1.6){$b_0=0$}
\rput(4.8,5){\scriptsize{value}}
\rput(4.8,4.7){\scriptsize{=$1/2$}}
\psline[linestyle=dashed, dash=1pt 1pt]{->}(4.8,4.54)(4.65,4.34)

\rput(7.2,3.5){\small0}
\pscircle(7.2,3.5){1.5}
\rput(7.2,1.6){$|b_0|<2/3$}
\pscustom[linewidth=0.7pt,fillstyle=hlines]{\psarc(7.2,3.5){0.35}{5}{45}
\psarc(7.2,3.5){0.6}{50}{0}}
\psarc(7.2,3.5){0.9}{-4}{54}
\pscurve(7.55,3.55)(7.6875,3.55)(7.825,3.5)(7.925,3.5)(8.1,3.42)

\pscurve(7.42,3.76)(7.5025,3.83)(7.605,3.985)(7.67,4.05)(7.73,4.2)
\psarc(7.2,3.5){1.15}{-15}{65}
\psarc(7.2,3.5){0.9}{65}{-15}
\pscurve(7.58,4.3)(7.65,4.42)(7.66,4.55)
\pscurve(8.06,3.29)(8.22,3.25)(8.32,3.17)
\pscircle(7.2,3.5){1.35}
\pscircle*(7.77,3.77){0.03}
\rput(8.72,4.48){\small$b_0$}
\psarc(7.2,3.5){0.74}{15}{35}
\psarc(7.2,3.5){0.53}{15}{35}
\pscurve(7.71,3.63)(7.76,3.66)(7.81,3.665)(7.86,3.7)(7.91,3.7)
\pscurve(7.65,3.81)(7.6875,3.82)(7.725,3.87)(7.7625,3.88)(7.8,3.93)
\pscurve[linestyle=dashed, dash=1pt 1pt]{->}(8.55,4.4)(8,4)(7.8,3.8)

\rput(11,3.5){\small0}
\rput(11,1.6){$|b_0|=1$}
\pscircle(11,3.5){1.5}
\pscustom[linewidth=0.7pt,fillstyle=hlines]{\psarc(11,3.5){0.75}{-55}{102}
\pscurve(10.85,4.22)(10.64,4.55)(10.4,4.84)
\psarc(11,3.5){1.5}{115}{-67}
\pscurve(11.58,2.16)(11.54,2.505)(11.41,2.89)}
\psarc(11,3.5){1}{-25}{75}
\pscurve(11.3,4.97)(11.32,4.65)(11.26,4.43)
\pscurve(11.9,3.11)(12.15,2.87)(12.2,2.67)
\psarc(11,3.5){1.2}{0}{50}
\pscurve(11.75,4.4)(11.83,4.52)(11.85,4.71)
\pscurve(12.2,3.49)(12.35,3.48)(12.5,3.41)
\rput(12.57,4.15){\small$b_0$}
\pscircle*(12.36,4.1){0.04}
\end{pspicture}
\caption{Level sets of the function $b\to d(b,b_0)$}
\begin{quote}{\small Each level curve of $b\to d(b,b_0)$ is a union of curves of the form $\{c_1\exp(\pi i x_2):c_2\le x_2\le c_3\}$ and $\{(c_1\pm x_2)\exp(\pi i x_2):c_2\le x_2\le c_3\}$ if the value is not $|b_0|/2$.}\end{quote}
\end{figure}
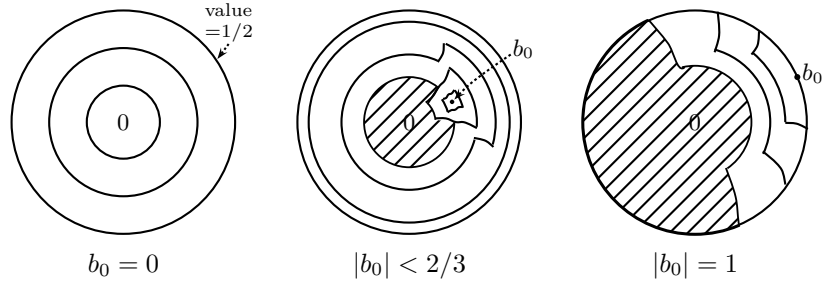

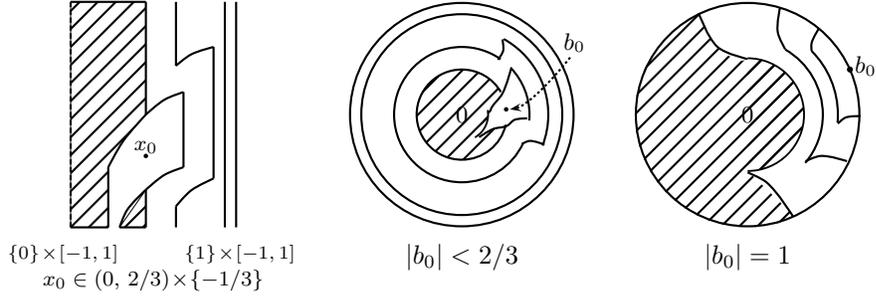
\begin{figure}[!h]
\begin{pspicture}(-0.85,1.3)(10.5,5)
\rput(3,3.07){\footnotesize$x_0$}
\pscircle*(3,2.95){0.03}
\psline[linestyle=dashed, dash=2pt 1pt](2,2)(2,5)
\pscustom[linewidth=0.1pt]{\psline[linestyle=dashed, dash=0.1pt 0.1pt](2,2)(2,5)\gsave
\psline(3,5)(3,3.5)
\pscurve(3,3.5)(2.7,3.15)(2.5,2.8)
\psline(2.5,2.8)(2.5,2)
\fill[linewidth=0.1pt,fillstyle=hlines]\grestore}
\psline(2,2)(2.5,2)
\psline(2,5)(3,5)
\psline(2.65,2)(3,2)
\psline(3,5)(3,3.5)
\pscurve(3,3.5)(2.7,3.15)(2.5,2.8)
\psline(2.5,2.8)(2.5,2)
\pscurve(2.65,2)(2.75,2.23)(3,2.5)
\pscustom[linewidth=0.1pt]{\psline[linestyle=dashed, dash=2pt 1pt](2.65,2)(3,2.5)\gsave
\psline(3,2.5)(3,2)\fill[linewidth=0.1pt,fillstyle=hlines]\grestore}
\psline(3,2.5)(3,2)
\pscurve(3,3.5)(3.25,3.7)(3.5,3.8)
\psline(3.5,3.8)(3.5,2.8)
\pscurve(3.5,2.8)(3.25,2.7)(3,2.5)
\psline(3.4,5)(3.4,4.15)
\pscurve(3.4,4.15)(3.6,4.3)(3.9,4.4)
\psline(3.9,4.4)(3.9,2.65)
\pscurve(3.9,2.65)(3.6,2.49)(3.4,2.28)
\psline(3.4,2.28)(3.4,2)
\psline(4.2,5)(4.2,2)
\psline(4.05,5)(4.05,2)
\rput(1.9,1.65){\scriptsize$\{0\}\!\times\![-1,1]$}
\rput(4.25,1.65){\scriptsize$\{1\}\!\times\![-1,1]$}
\rput(3.1,1.3){\footnotesize$x_0\in(0,\,2/3)\!\times\!\{-1/3\}$}

\rput(7.2,3.5){\small0}
\pscircle(7.2,3.5){1.5}
\rput(7.2,1.6){$|b_0|<2/3$}
\pscustom[linewidth=0.7pt,fillstyle=hlines]{\psarc(7.2,3.5){0.35}{-60}{10}\gsave
\psarc(7.2,3.5){0.6}{30}{-13}\fill[linewidth=0.1pt,fillstyle=hlines]\grestore}
\psarc(7.2,3.5){0.6}{30}{-13}
\psarc(7.2,3.5){0.9}{-5}{43}
\pscurve(7.43,3.2)(7.62,3.31)(7.8,3.35)(7.9,3.43)(8.07,3.45)
\pscurve(7.54,3.52)(7.57,3.54)(7.6,3.6)(7.74,3.78)(7.86,4.13)
\psarc(7.2,3.5){1.15}{-21}{61}
\psarc(7.2,3.5){0.9}{55}{-25}
\pscurve(7.68,4.23)(7.75,4.35)(7.74,4.5)
\pscurve(8,3.11)(8.1,3.15)(8.3,3.1)
\pscircle(7.2,3.5){1.35}
\rput(8.7,4.46){\small$b_0$}
\pscircle*(7.79,3.57){0.03}
\pscurve[linestyle=dashed, dash=1pt 1pt]{->}(8.64,4.26)(8,3.62)(7.84,3.57)

\rput(11,3.5){\small0}
\rput(11,1.6){$|b_0|=1$}
\pscircle(11,3.5){1.5}
\pscustom[linewidth=0.7pt,fillstyle=hlines]{
\psarc(11,3.5){0.75}{-90}{90}\gsave
\pscurve(11,4.24)(10.6,4.4)(10.35,4.85)
\psarc(11,3.5){1.5}{115}{-65}\fill[linewidth=0.1pt,fillstyle=hlines]\grestore}
\pscurve(11,4.24)(10.6,4.4)(10.35,4.85)
\pscurve(11,2.73)(11.25,2.6)(11.6,2.15)
\psarc(11,3.5){1}{-40}{60}
\pscurve(11.43,4.92)(11.4,4.6)(11.5,4.37)
\pscurve(11.77,2.88)(11.97,2.95)(12.32,2.88)
\psarc(11,3.5){1.2}{-5}{45}
\pscurve(11.85,4.35)(11.85,4.45)(11.95,4.65)
\pscurve(12.2,3.38)(12.3,3.46)(12.5,3.48)
\rput(12.57,4.15){\small$b_0$}
\pscircle*(12.36,4.1){0.04}
\end{pspicture}
\caption{Level sets of $x\to\tilde d_X(x,x_0)$ and $b\to\tilde d(b,b_0)$}
\begin{quote}{\small Each level curve of $b\to\tilde d(b,b_0)$ is a union of curves of the form $\{c_1\exp(\pi i x_2):c_2\le x_2\le c_3\}$ and $\{x_1\exp\big(\pi i(c_4\pm x_1+\phi(x_1))\big):c_5\le x_2\le c_6\}$ if the value is not $|b_0|/2$.}\end{quote}
\end{figure}
Under these conditions, we will prove that there are only two isometries from $\tilde B$ to $B$ which carries $(1,0)$ to $b_0$, they are both non-subanalytic, one is orientation-preserving and coincides with $\Phi_{b_0}$, and the other is orientation-reversing.
We consider only the orientation-preserving one because $\rho\circ\Phi'$ is the orientation-preserving one by the uniform property of $d_X$ for an orientation-reversing one $\Phi'$ and the $\rho\in O(2)$ such that $\rho(b_0)=b_0$ and the determinant of $\rho$ is $-1$.
First $\Phi_{b_0}$ is non-subanalytic for the following reason:\par
As above, let $\Phi_X:\tilde X\to X$ be the isometry which induces $\Phi_{b_0}$, and $x_0\in\R$ such that $f(1,x_0)=b_0$.
Then $\exp(\pi ix_0)=b_0$ and
$$\Phi_X((0,\,1]\times\{0\})=\{(x_1,x_0+\phi(x_1)):x_1\in(0,\,1]\}.$$\noindent
Hence\vskip-10mm
\begin{align*}\Phi_{b_0}(b_0*0)-\{0\}=f\{(x_1,x_0+\phi(x_1)):x_1\in(0,\,1]\}=\qquad\\
\{x_1\exp(\pi i(x_0+\phi(x_1))):x_1\in(0,\,1]\}=\{b_0x_1\exp(\pi i\phi(x_1)):x_1\in(0,\,1]\},\end{align*}
where $b_0*0$ denotes the segment with ends $b_0$ and 0 in $B$.
If $x_1$ moves from 1 to 0 then $b_0x_1\exp(\pi i\phi(x_1))$ spins infinitely around 0 and converges to 0.
Thus the set $\Phi_{b_0}(b_0*0)$ and the map $\Phi_{b_0}$ are non-subanalytic.\par
Lastly, we see that if an orientation-preserving isotopy $\Phi:\tilde B\to B$ carries $(1,0)$ to $b_0$ then $\Phi=\Phi_{b_0}$.
For this, it suffices to show the following statement by the uniform property of $d_X$:\vskip1mm\noindent
{\it Let $\Psi:(B,d)\to(B,d)$ be an orientation-preserving isotopy such that $\Psi(b_0)=b_0$.
Then $\Psi=\id$.}\par
Moreover, we only need to see that if $\Psi=\id$ at some $b\in\partial B$ then it holds on a neighborhood of $b$ in $\partial B$ and if $\Psi=\id$ on $b*b'$ for some $b'\in b*0-\{0\}$ then it holds on $b*b''$ for some $b''\in b'*0-\{b'\}$.
We can translate these statements to the following ones on $d_X$:\vskip1mm\noindent
{\it Given an orientation-preserving isometry $\Psi_X:(X,d_X)\to(X,d_X)$, if $\Psi_X=\id$ at some $(1,x_2)\in\{1\}\times\R$ then it holds on a neighborhood of $(1,x_2)$ in $\{1\}\times\R$, and if $\Psi_X=\id$ on $[x_1,\,1]\times\{x_2\}$ for some $(x_1,x_2)\in(0,\,1]\times\R$ then it holds on $[x'_1,\,1]\times\{x_2\}$ for some $x'_1\in(0,\,x_1)$.}\par
To prove the first statement, we take a number $x'_2$ in $\R$ larger than $x_2$ and near $x_2$.
Then $(1,x'_2)$ is uniquely determined as an element of $\{1\}\times\R$ by its Euclidean distance to $(1,x_2)$ by the following equality:
$$d_X(x,x')=\max\{|x_1-x'_1|,|x_2-x'_2|\}\ \text{ for }x=(x_1,x_2)\in X\text{ and }x'=(x'_1,x'_2)\text{ close to }x.$$
Hence $\Psi_X(1,x'_2)=(1,x'_2)$ since $\Psi_X(1,x_2)=(1,x_2)$ and $\Psi_X$ is metric preserving.
We obtain the same equality for $x'_2$ smaller than $x_2$.
For the second, we take a number $x''_1$ in $(x_1,\,1)$ near $x_1$.
Then any number in $[2x_1-x''_1,\,x''_1]\times\{x_2\}$ is determined by its Euclidean distances to $(x_1,x_2)$ and $(x''_1,x_2)$ by the same reason, and we have $\Psi_X=\id$ at the point since $\Psi_X=\id$ at $(x_1,x_2)$ and $(x''_1,x_2)$.
Thus the statements are proved, and we see that $B$ and $\tilde B$ are isometric subanalytic metric spaces which are not subanalytically isometric.\par
It only remains to modify $d$ and $\tilde d$ to semialgebraic metrics.
As $d_X$ and $\tilde d_X$ are semialgebraic, we need to modify the $f$.
Let $g:[0,\,2]\to\partial B$ be a semialgebraic $C^1$ approximation of $f|_{\{1\}\times[0,\,2]}$ the map $[0,\,2]\ni x_2\to\exp(\pi i x_2)\in\partial B$, such that $g(0)=g(2)$.
Set $f'(x_1,x_2)=x_1g(x_2-2n)$ for $(x_1,x_2)\in(0,\,1]\times[2n,\,2n+2],\ n\in\Z$.
Then $f'$ is semialgebraic on $(0,\,1]\times[2n,\,2n+2]$ for any $n$, it is a covering map, and Conditions $(*),\,(**),\,\widetilde{(*)}$ and $\widetilde{(**)}$ hold.
Hence $d'$ and $\tilde d'$, defined by $d_X,\,\tilde d_X$ and $f'$ likewise $d$ and $\tilde d$ by $d_X,\,\tilde d_X$ and $f$, are semialgebraic metrics on $B$.
We see that $(B,d')$ and $(B,\tilde d')$ are isometric but not semialgebraically isometric as in the subanalytic case.
Thus Theorem 2 is proved.\qed
\vskip2mm\noindent
{\it Proof of Theorem 3 in the case $r>0$.} To avoid confusion, we write $\psi:Y\to X$ in Theorem 3 as $\phi:B\to A$.
Assume that $r=1$ and the set of definably $C^1$ smooth points of $A$ of dimension $k$ is dense in $A$ for simplicity of notation.
Let $\Reg A$ denote this set.
Embed $A$ in $R^n$ so that $A$ is bounded in $R^n$.
First we choose a definable $C^1$ manifold possibly with corners as $B$, and later we smooth the corners if exists.
For this, it suffices to prove the following statement:\vskip1mm\noindent
Statement 1. {\it Let $X$ be a compact definable set, and $X_1$ a definable closed nowhere dense subset of $X$.
Then there exist a compact definable $C^1$ manifold $Y$ possibly with corners, a definable subset $Y_1$ of $Y$ and a definable $C^0$ map $\psi:Y\to X$ such that $\psi|_{\Int Y}$ is a $C^1$ diffeomorphism onto $\Reg(X-X_1)$, $\psi(Y_1)=X_1$, $\psi(\partial Y)=X-\Reg(X-X_1)$ and $Y_1$ is the union of the interiors of some faces of $Y$ of dimension $k-1$.}\par
Indeed, assume that Statement 1 holds, and apply it to $\overline A$ and $(\overline A-A\overline)$ for $A$ in Theorem 3 as $X$ and $X_1$ in Statement 1 respectively.
Let $Y,\,Y_1$ and $\psi:Y\to X$ be consequences of Statement 1.
Then $\psi^{-1}((\overline A-A\overline))$ is a closed definable subset of $Y$, hence $Y-\psi^{-1}((\overline A-A\overline))$ is a definable $C^1$ manifold possibly with corners, and $\big(Y-\psi^{-1}((\overline A-A\overline))\big)\cup Y_1$ is also so because $Y_1\cup\psi^{-1}((\overline A-A\overline))\subset\partial Y$.
Furthermore, $\psi\big(\Int((Y-\psi^{-1}((\overline A-A\overline)))\cup Y_1)\big)=\Reg A$ because
$$\Int((Y-\psi^{-1}((\overline A-A\overline)))\cup Y_1)=\Int(Y-\psi^{-1}((\overline A-A\overline)))=\Int Y\quad\text{and}$$
$$\psi(\Int Y)\!=\!\Reg(X\!-\!X_1)\!=\!\Reg(\overline A\!-\!(\overline A\!-\!A\overline))\!=\!\Reg\!\big(A\!-\!((\overline A\!-\!A\overline)\!-\!(\overline A\!-\!A))\big)\!=\!\Reg A.$$
Thus the set $\big(Y-\psi^{-1}((\overline A-A\overline))\big)\cup Y_1$ and the restriction of $\psi$ to it satisfy the conditions on $B$ and $\phi$, respectively, in Theorem 3.\par
Now we prove Statement 1.
Let $\{C_i:i\}$ be a Whitney definable $C^1$ stratification of $X$ such that $\dim C_i=i$, $C_k=\Reg(X-X_1)$ and $X_1$ is the union of some $C_i$'s.
Let $\{U_i=(|U_i|,\pi_i,\xi_i):i\}$ be a controlled definable tube system for $\{C_i\}$, and $\{\epsilon_i:i\}$ the set of sufficiently small positive numbers such that $\epsilon_i\gg\epsilon_{i'}$ if $i<i'$.
For each $i=0,\ldots,k-1$, set
$$V_i=\{x\in U_i:\xi_i(x)\le\epsilon_i\}-\cup_{i'<i}\{x\in U_{i'}:\xi_{i'}(x)<\epsilon_{i'}\}.$$
Then the restrictions $\pi_i|_{V_i}:V_i\to(C_i-\cup_{i'<i}V_{i'}\overline)$ and $\xi_i|_{V_i}$ are a definable retraction and a definable continuous function, respectively, and they induce a hole-blow-up $\psi_i:(X-\cup_{i'\le i}V_{i'}\overline)\to(X-\cup_{i'<i}V_{i'}\overline)$.
Here $\psi_i$ is originally defined on $(\{x\in X\cap U_i:\epsilon_i\le\xi_i(x)\le2\epsilon_i\}-\cup_{i'\le i}V_{i'}\overline)$ so that $\psi_i=\id$ on $(\{x\in X\cap U_i:\xi_i(x)=2\epsilon_i\}-\cup_{i'\le i}V_{i'}\overline)$, and hence we extend $\psi_i$ to $(X-\cup_{i'\le i}V_{i'}\overline)$ by setting $\psi_i=\id$ outside the original domain.
Thus we have a sequence of hole-blows-up $(X-\cup_{i<k}V_i\overline)\overset{\psi_{k-1}}{\longrightarrow}\cdots\overset{\psi_1}{\longrightarrow}(X-V_0\overline)\overset{\psi_0}{\longrightarrow}X$.
Set $Y=(X-\cup_{i<k}V_i\overline)$ and $\psi=\psi_0\circ\cdots\circ\psi_{k-1}$.
Then $Y$ is a compact definable $C^1$ manifold possibly with corners, $\psi$ is a definable $C^0$ map, and for each $i<k$, the set $\psi^{-1}_i(C_i)-\cup_{i'\not=i}V_{i'}$ is a definable $C^1$ manifold and the union of the interiors of some faces of $Y$ of dimension $k-1$.
Moreover, we have
\begin{align*}\psi(\psi^{-1}_i(C_i)-\cup_{i'\not=i}V_{i'})=\psi_0\circ\cdots\circ\psi_i(\psi^{-1}_i(C_i)-\cup_{i'<i}V_{i'}-\cup_{i<i'<k}C_{i'})\\
=\psi_0\circ\cdots\circ\psi_{i-1}(C_i-\cup_{i'<i}V_{i'})=C_i.\qquad\qquad\qquad\end{align*}
Hence if we let $Y_1$ be the union of the manifolds $\psi^{-1}_i(C_i)-\cup_{i'\not=i}V_{i'}$ such that $C_i\subset X_1$ then $\psi(Y_1)=X_1$, $\psi(\partial Y)=X-\Reg(X-X_1)$ and $Y_1$ is the union of the interiors of some faces of dimension $k-1$.
Thus Statement 1 is proved, and we obtain a definable definable $C^1$ manifold $Y$ possibly with corners which satisfies the conditions in Theorem 3.\par
It remains to smooth the corners of the manifold $Y$ possibly with corners.
Assume that there are the corners.
The smoothing is well known in the case $\R$.
That proof uses a vector field and its integration.
However, a vector field does not necessarily induce its integration in the general case.
Hence we directly construct a flow.
Let $f$ be a non-negative definable continuous function on $Y$ such that $f^{-1}(0)=\partial Y$, if the germ of $Y$ at a point of $\partial Y$ is definably diffeomorphic to the germ of $[0,\,\infty)^l\times R^{k-l}$ at 0 then $f(x_1,\ldots,x_k)$ is of the form $x_1\cdots x_l$ through some definable $C^1$ diffeomorphism and hence $f$ is $C^1$ regular on the intersection of $\Int Y$ and some neighborhood of $\partial Y$ in $Y$.
Let $\tilde Y$ be a definable $C^1$ manifold containing $Y$ and of the same dimension.
Then it suffices to find a definable $C^1$ 1-dimensional foliation $\{F_a:a\in A\}$ of a neighborhood of $\partial Y$ in $\tilde Y$ such that each leaf $F_a$ intersects transversally each face of $Y$ of dimension $k-1$ and the hypersurface $f^{-1}(t)$ for each $t\in(0,\,1]$ close to 0 for the following reason:\par
We need to define a definable $C^1$ manifold $Z$ with boundary and a definable homeomorphism $\eta:Z\to Y$ such that $\eta|_{\Int Z}$ is a $C^1$ diffeomorphism onto $\Int Y$.
Shrink the domain where the foliation is defined, and choose $t_0\in(0,\,1]$ close to 0 so that $f^{-1}((0,\,t_0))$ is the domain, the restriction $f|_{f^{-1}((0,\,t_0))}:f^{-1}((0,\,t_0))\to(0,\,t_0)$ is $C^1$ regular and proper and each $F_a$ intersects transversally $f^{-1}(t)$ for $t\in(0,\,t_0]$.
Let $\chi:[t_0/2,\,t_0]\to[0,\,t_0]$ be a definable $C^1$ diffeomorphism such that $\chi(t_0/2)=0$ and $d\chi/d t=1$ on $[2t_0/3,\,t_0]$.
Set $Z=f^{-1}([t_0/2,\,\infty))$ and define a definable homeomorphism $\eta:Z\to Y$ so that $\eta=\id$ on $f^{-1}([t_0,\,\infty))$ and for $z\in f^{-1}((0,\,t_0])$, $\eta(z)$ is the point in $F_a$ such that $f\circ\eta(z)=\chi\circ f(z)$, where $F_a$ is such that $z\in F_a$, i.e., $z$ moves to $\eta(z)$ along $F_a$ so that the value of $f$ moves from $t$ to $\chi(t)$.
Then $\eta$ is well defined, and $\eta|_{\Int Z}$ is a $C^1$ diffeomorphism onto $\Int Y$, which explains the reason.\par
We construct $\{F_a\}$ by induction.
For simplicity of notation and without loss of generality, we assume that the face assumption holds for $Y$ and each proper face of $Y$ is contractible.
Let $\{D_i:i=0,\ldots,k-1\}$ denote the canonical stratification of $\partial Y$, i.e., $\{\Int Y\}\cup\{D_i:i=0,\ldots,k-1\}$ is the canonical stratification of $Y$.
Consider the basis step.
We can regard the triple of a point $y$ of $D_0$ and its definable neighborhoods in $Y$ and $\tilde Y$ as the one of 0, $[0,\,\infty)^k$ and $R^k$ through a definable $C^1$ diffeomorphism.
Hence we can translate the problem of the construction of $\{F_a\}$ around $y$ to the one around 0 in $R^k$.
We define a foliation of the set $(-1,\,1)^k$ to be the family of the lines parallel to the line passing through 0 and $(1,\ldots,1)$.
Translate this foliation back to a foliation of a neighborhood of $y$ in $\tilde Y$ through the above definable $C^1$ diffeomorphism.
Thus we obtain a definable neighborhood $W_0$ of $D_0$ in $\tilde Y$ and a definable $C^1$ 1-dimensional foliation $\{F_a:a\in A_0\}$ of $W_0$.
Clearly this foliation satisfies the transversality condition.
Assume as the induction hypothesis that there are a definable neighborhood $W_{l-1}$ of $D_0\cup\cdots\cup D_{l-1}$ and a definable $C^1$ 1-dimensional foliation $\{F_a:a\in A_{l-1}\}$ of $W_{l-1}$ such that the transversality condition hold.\par
We will shrink $W_{l-1}$ and extend $W_{l-1}$ and $\{F_a:a\in A_{l-1}\}$ to $W_l$ and $\{F_a:a\in A_l\}$, respectively, i.e., we modify $\{F_a:a\in A_{l-1}\}$ outside a smaller neighborhood of $D_0\cup\cdots\cup D_{l-1}$ and extend it to $\{F_a:a\in A_l\}$.
Set $D'_l=\overline{D_l}\cap W_{l-1}$.
Here we can assume that $D_l$ is connected, and we regard a definable neighborhood of $D_l$ in $Y$ as $D_l\times[0,\,\infty)^{k-l}$ since $D_l$ is contractible.
Then the problem is reduced to the following statement:\vskip1mm\noindent
Statement 2. {\it Let $D'_l$ be a definable neighborhood of $\overline{D_l}\cap(D_0\cup\cdots\cup D_{l-1})$ in $\overline{D_l}$.
Let $f$ be the function on $\overline{D_l}\times R^{k-l}$ defined so that $f(y,x_1,\ldots,x_{k-l})=x_1\cdots x_{k-l}$ for $(y,x_1,\ldots,x_{k-l})\in\overline{D_l}\times R^{k-l}$.
Let $\{F_a:a\in A_{l-1}\}$ be a definable $C^1$ 1-dimensional foliation of a definable neighborhood $W_{l-1}$ of $D'_l\times\{0\}$ in $D'_l\times R^{k-l}$ such that $F_a\subset\{0\}\cup D'_l\times((0,\,\infty)^{k-l}\cup(-\infty,\,0)^{k-l})$ if $F_a$ passes through $D'_l\times\{0\}$ and each $F_a$ intersects transversally $f^{-1}(t)$ for each $t\not=0$ and $D'_l\times R^{l'}\times\{0\}\times R^{k-l-l'-1},\ l'=0,\ldots,k-l-1$.
Then shrinking $D'_l$ and $W_{l-1}$, we can extend $\{F_a:a\in A_{l-1}\}$ to a definable $C^1$ 1-dimensional foliation $\{F_a:a\in A_l\}$ of some definable neighborhood $W_l$ of $\overline{D_l}\times\{0\}$ in $\overline{D_l}\times R^{k-l}$ so that the same conditions continue to hold.}\par
We naturally define a definable $C^1$ 1-dimensional foliation $\{F_b:b\in B\}$ of $D_l\times N$ so that for each $x\in D_l$, $\{F_b\cap\{x\}\times R^{k-l}:b\in B\}$ is the product of $\{x\}$ and the foliation of $N$ defined in the same way as in the basis step, where $N$ is some neighborhood of 0 in $R^{k-l}$.
Set $W_l=D_l\times N$.
Then shrinking $W_{l-1},\ D'_l$ and $W_l$ if necessary, we obtain a definable homeomorphism $\tau:D'_l\times N\to W_{l-1}\cap D_l\times R^{k-l}$ such that $\tau=\id$ on $D'_l\times(N\cap\cup_{l'=0}^{k-l-1}[0,\,\infty)^{l'}\times\{0\}\times[0,\,\infty)^{k-l-l'-1})$, $f\circ\tau=f$ and each $F_b$ in $D'_l\times N$ is carried to the $F_a$ for $a\in A_{l-1}$ such that $F_b\cap D'_l\times[0,\,\infty)^{k-l}=F_a\cap D'_l\times[0,\,\infty)^{k-l}$.
The $\tau$ has the following property of differentiability:\par
Let $L$ denote the line in $R^{k-l}$ passing through 0 and $(1,\ldots,1)$, and $M_{l'}$ the union of the lines passing through $[0,\,\infty)^{l'}\times\{0\}\times[0,\,\infty)^{k-l-l'-1}$ and parallel to $L$ for $l'=0,\ldots,k-l-1$.
For example, if $k-l=2$,
$$M_0=\{(x_1,x_2)\in R^2:x_1\le x_2\}\ \text{ and }\ M_1=\{(x_1,x_2)\in R^2:x_1\ge x_2\}.$$
Then each $M_{l'}$ is linearly homeomorphic to $[0,\,\infty)^{k-l-1}\times R$, and the family of the strata of the canonical stratifications of $M_{l'}$ for all $l'$ is a stratification of $R^{k-l}$ into open subsets of linear spaces.
The property is that $\tau$ is a piecewise $C^1$ diffeomorphic with respect to $\{D'_l\times M_{l'}:l'\}$, i.e., $\tau|_{D'_l\times M_{l'}}$ is a $C^1$ diffeomorphism onto its image for any $l'$.\par
It suffices to extend $\tau$ to a definable homeomorphism $\tilde\tau:D_l\times N\to(W_{l-1}\cap D_l\times R^l)\cup(D_l\times N-W_{l-1})$ so that $\tilde\tau=\id$ on $D_l\times(N\cap\cup_{l'=0}^{k-l-1}[0,\,\infty)^{l'}\times\{0\}\times[0,\,\infty)^{k-l-l'-1})$ and $\tilde\tau$ is a piecewise $C^1$ diffeomorphic with respect to $\{D_l\times M_{l'}:l'\}$.
We clearly set $\tilde\tau=\id$ on $D_l\times(N\cap\cup_{l'=0}^{k-l-1}[0,\,\infty)^{l'}\times\{0\}\times[0,\,\infty)^{k-l-l'-1})$ and outside a small definable neighborhood of $D'_l\times N$ in $D_l\times N$.
Then by shrinking $D'_l$ and using a definable $C^1$ partition of unity, we obtain the $\tilde\tau$ because of the next obvious fact:
Let $\tau_1$ be a definable orientation-preserving $C^1$ diffeomorphism of $D_l\times[0,\,\infty)^{k-l-1}\times R$ such that $\tau_1=\id$ on $D_l\times[0,\,\infty)^{k-l-1}\times\{0\}$, and let $\theta$ be a definable $C^1$ function on $D_l$ such that $0\le\theta\le1$.
Let $\id$ denote the identity map of $D_l\times[0,\,\infty)^{k-l-1}\times R$.
Then the fact is that the map $D_l\times[0,\,\infty)^{k-l-1}\times R\ni(x,y,t)\to\theta(x)\tau_1(x,y,t)+(1-\theta(x))\id(x,y,t)\in D_l\times[0,\,\infty)^{k-l-1}\times R$ restricted to some definable neighborhood of $D_l\times[0,\,\infty)^{k-l-1}\times\{0\}$ is a $C^1$ embedding.
Hence Statement 2 and the smoothing are proved, which completes the proof of Theorem 3 in the case $r>0$.\qed
\vskip2mm\noindent
{\it Proof of Theorem 3 in the case $r=0$.} Let $\Reg X$ denote the set of definably $C^0$ smooth points of $X$ of dimension $k$.
Assume that $X$ is bounded in $R^n$ and $\Reg X$ is dense in $X$ for simplicity of notation.
Moreover, we suppose that there is a finite simplicial complex $K$ such that its underlying polyhedron is contained in $R^n$ and $X$ is the union of the interiors of some simplices of $K$ by the triangulation theorem of definable sets.
Then $\Reg X$ is the subset of $X$ of PL smooth points of $X$ of dimension $k$, i.e., the subset of $X$ of points where the germ of $X$ is PL homeomorphic to the one of $R^k$ at 0 by Theorem 1.1 (o-minimal Hauptvermutung) in \cite{S2}.
Note that $\{\Int\sigma:\sigma\in K\}$ is a Whitney definable $C^1$ stratification of $\overline X$.
Set $\{C_i:i\}=\{\Int\sigma:\sigma\in K,\Int\sigma\subset X-\Reg X\}$, and gather $C_i$ of the same dimensions so that $\dim C_i=i$.
We define a controlled definable tube system $\{U_i=(|U_i|,\pi_i,\xi_i):i\}$ for $\{C_i\}$, a set of sufficiently small positive numbers $\{\epsilon_i:i\}$ such that $\epsilon_i\gg\epsilon_{i'}$ if $i<i'$ and a set of subsets of $R^n$ $\{V_i:i\}$ as in the case $r>0$.
Then we have a sequence of hole-blows-up $(X-\cup_{i<k}V_i\overline)\overset{\psi_{k-1}}{\longrightarrow}\cdots\overset{\psi_1}{\longrightarrow}(X-V_0\overline)\overset{\psi_0}{\longrightarrow}\overline X$.
Set $Y=\psi^{-1}(X)\cap(X-\cup_{i<k}V_i\overline)$ and $\psi=\psi_0\circ\cdots\circ\psi_{k-1}|_Y$.
Then $Y$ and $\psi$ satisfy the conditions in Theorem 3 by the same reason as in the case $r>0$.
Thus we prove Theorem 3 in the case $r=0$.\qed

\small\quad\ \ Graduate School of Mathematics, Nagoya University, Chikusa, Nagoya 464-8602, Japan\vskip1mm
{\it E-mail address}: shiota@math.nagoya-u.ac.jp.
\end{document}